\documentclass[10pt]{article}
\usepackage{amsmath,amsfonts,amssymb,amsthm,amsxtra,latexsym,verbatim}
\usepackage{bookman}
\usepackage[mathscr]{eucal}
\usepackage[all]{xy}

\theoremstyle{plain}
\newtheorem{thm}{Theorem}[section]
\newtheorem{cor}[thm]{Corollary}
\newtheorem{lemma}[thm]{Lemma}
\newtheorem{prop}[thm]{Proposition}

\theoremstyle{definition}
\newtheorem{ass}[thm]{Assertion}

\newtheorem{defn}[thm]{Definition}

\theoremstyle{remark}

\newtheorem{xmpl}[thm]{Example}

\newtheorem{rmk}[thm]{Remark}

\newcommand{\C}{\mathbb{C}}
\newcommand{\N}{\mathbb{N}}
\newcommand{\Q}{\mathbb{Q}}
\newcommand{\R}{\mathbb{R}}
\newcommand{\Z}{\mathbb{Z}}

\renewcommand{\P}{\mathbf{P}}

\newcommand{\A}{\mathcal{A}}
\newcommand{\E}{\mathcal{E}}
\newcommand{\F}{\mathcal{F}}

\newcommand{\J}{\mathcal{J}}
\newcommand{\K}{\mathcal{K}}
\renewcommand{\L}{\mathcal{L}}
\newcommand{\M}{\mathcal{M}}
\renewcommand{\O}{\mathcal{O}}

\renewcommand{\a}{\mathfrak{a}}
\newcommand{\m}{\mathfrak{m}}
\newcommand{\p}{\mathfrak{p}}

  \DeclareMathOperator{\gr}{gr} 
  
 \DeclareMathOperator{\rad}{rad} \DeclareMathOperator{\rank}{rank}
\DeclareMathOperator{\Sh}{Sh} \DeclareMathOperator{\Proj}{Proj} \DeclareMathOperator{\Spec}{Spec}
\DeclareMathOperator{\Supp}{Supp}  \DeclareMathOperator{\Tor}{Tor}

\begin{document}

\title{Toric singularities revisited}
\author{Howard~M Thompson}

\maketitle
\begin{abstract}
   In \cite{kK94}, Kato defined his notion of a log regular scheme and studied the local behavior of such schemes. A
   toric variety equipped with its canonical logarithmic structure is log regular. And, these schemes allow one to
   generalize toric geometry to a theory that does not require a base field. This paper will extend this theory by
   removing normality requirements.
\end{abstract}


\tableofcontents

\section*{Conventions and Notation}

All monoids considered in this paper are commutative and cancellative. All rings considered in this paper are
commutative and unital. See Kato~\cite{kK89} for an introduction to log schemes. There Kato defines pre-log structures
and log structures on the \'etale site of $X$. However, we will use the Zariski topology throughout this paper.

\begin{list}{}{\setlength{\itemindent}{0pt}\setlength{\leftmargin}{0.8in}\setlength{\labelwidth}{0.8in}}
   \item[$P^*$ \hfill] the unit group of the monoid $P$.
   \item[$\overline{P}$ \hfill] the sharp image of the monoid $P$, $\overline{P}=P/P^*$ is the orbit space under
      the natural action of $P^*$ on $P$.
   \item[$P^+$ \hfill] $P^+=P\setminus P^*$ is the maximal ideal of the monoid $P$.
   \item[$P^{gp}$ \hfill] the group generated by $P$, that is, image of $P$ under the left adjoint of the inclusion
      functor from Abelian groups to monoids.
   \item[$P^{sat}$ \hfill] the saturation of $P$, that is, $\{p\in P^{gp}\mid np\in P\text{ for some }n\in\N^+\}$.
   \item[${R[P]}$ \hfill] the monoid algebra of $P$ over a ring $R$.  The elements of $R[P]$ are written as
      ``polynomials''. That is, they are finite sums $\sum r_pt^p$ with coefficients in $R$ and exponents in $P$.
   \item[$(K)$ \hfill] the ideal $\beta(K)A$, where $\beta:P\to A$ is a monoid homomorphism with respect to
      multiplication on $A$ and $K$ is an ideal of $P$. We say such an ideal is a \emph{log ideal} of $A$.
   \item[${R[[P]]}$ \hfill] the $(P^+)$-adic completion of $R[P]$.
   \item[$P_{\p}$ \hfill] the localization of the monoid $P$ at the prime ideal $\p\subseteq P$.
   \item[$\dim P$ \hfill] the (Krull) dimension of the monoid $P$.
\end{list}

\section*{Introduction}\label{Ch:intro}\addcontentsline{toc}{section}{Introduction}

A toric variety is a normal irreducible separated scheme $X$, locally of finite type over a field $k$, which
contains an algebraic torus $T\cong(k^*)^d$ as an open set and is endowed with an algebraic action $T\times X\to
X$ extending the group multiplication $T\times T\to T$. According to Oda~\cite{tO91}:

\begin{quotation}
   \noindent The theory was started at the beginning of 1970's by Demazure ~\cite{mD70} in connection with
   algebraic subgroups of the Cremona groups, by Mumford et al.\ \cite{KKMS73} and Satake~\cite{iS73} in
   connection with compactifications of locally symmetric varieties, and by Miyake and Oda~\cite{MO75}. We were
   inspired by Hochster~\cite{mH72} as well as Sumihiro~\cite{hS74,hS75}.

   Comprehensive surveys from various different perspectives can be found in Danilov~\cite{vD78}, Mumford et al.
   \cite{KKMS73}, ~\cite{AMRT75} as well as~\cite{tO78,tO88}.
\end{quotation}

In~\cite{kK94}, Kato extended the theory of toric geometry over a field to an absolute theory, without base. This
is achieved by replacing the notion of a toroidal embedding introduced in \cite{KKMS73} with the notion of a log
structure. A toroidal embedding is a pair $(X,U)$ consisting of a scheme $X$ locally of finite type and an open
subscheme $U\subset X$ such that $(X,U)$ is isomorphic, locally in the \'etale topology, to a pair consisting of
a toric variety and its algebraic torus. Toroidal embeddings are particularly nice locally Noetherian schemes with
distinguished log structures.

A log structure on a scheme $X$, in the sense of Fontaine and Illusie, is a morphism of sheaves of monoids
$\alpha:\M_X\to\O_X$ restricting to an isomorphism $\alpha^{-1}(\O_X^*)\cong\O_X^*$. The theory of log structures
on schemes is developed by Kato in~\cite{kK89}. Log structures were developed to give a unified treatment of the
various constructions of deRham complexes with logarithmic poles. In~\cite{lI94} Illusie recalls the question
that motivated their definition:

\begin{quotation}
   \noindent Let me briefly recall what the main motivating question was. Suppose $S$ is the spectrum of a
   complete discrete valuation ring $A$, with closed (resp. generic) point $s$ (resp. $\eta$), and $X/S$ is a
   scheme with semi-stable reduction, which means that, locally for the \'etale topology, $X$ is isomorphic to
   the closed subscheme of $\mathbb{A}^n_S$ defined by the equation $x_1\cdots x_n=t$, where $x_1,\dots,x_n$ are
   coordinates on $\mathbb{A}^n$ and $t$ is a uniformizing parameter of $A$. Then $X$ is regular, $X_\eta$ is
   smooth, and $Y=X_s$ is a divisor with normal crossings on $X$. In this situation, one can consider, with
   Hyodo, the relative deRham complex of $X$ over $S$ with logarithmic poles along $Y$,
   \[
      \omega_{X/S}\spdot=\Omega_{X/S}\spdot(\log Y)
   \]
   (\cite{oH88} see also~\cite{oH91,HK94}). Its restriction to the generic fiber is the usual deRham complex
   $\Omega_{X_\eta/\eta}\spdot$ and it induces on $Y$ a complex
   \[
      \omega_Y\spdot=\O_Y\otimes\Omega_{X/S}\spdot(\log Y).
   \]
   One has $\omega_Y^i=\bigwedge^i\omega^1_Y$, and when $X$ is defined as above by $x_1\cdots x_r=t$,
   $\omega_Y^1$ is generated as an $\O_Y$-module by the images of the $d\log x_i\,(1\leq i\leq r)$ and
   $dx_i\,(r+1\leq i\leq n)$ subject to the single relation $\sum_{1\leq i\leq r}d\log x_i=0$. The analogue of
   $\omega_Y\spdot$ over $\C$ ($S$ is replaced by a disc) is the complex studied by Steenbrink in~\cite{jS76},
   which ``calculates'' $\mathbf{R}\psi(\C)$. If $X/S$ is of relative dimension $1$, $\omega^1_Y$ is simply the
   dualizing sheaf of $Y$ (which probably explains the notation $\omega$, chosen by Hyodo), and of course, in
   this case, $\omega_Y\spdot$ depends only on $Y$. In general, however, $\omega\spdot_Y$ depends not only on $Y$
   but also ``a little bit'' on $X$, so it is natural to ask: which extra structure on $Y$ is needed to define
   $\omega_Y\spdot$?

   Assume for simplicity that $Y$ has simple normal crossings, i.e. $Y$ is a sum of smooth divisors $Y_i$ meeting
   transversely. Let $\L_i$ be the invertible sheaf $\O_Y\otimes\O_X(-Y_i)$, and $s_i:\L_i\to\O_Y$ the map
   deduced by extension of scalars from the inclusion $\O_X(-Y_i)\hookrightarrow\O_X$. Then it is easily seen
   that the data consisting of the pairs $(\L_i,s_i)$, together with the isomorphism
   $\O_Y\xrightarrow{\sim}\bigotimes\L_i$ (coming from $Y=\sum Y_i$), suffice to define $\omega^1_Y$. Indeed,
   $\omega^1_Y$ can be defined as the $\O_Y$-module generated locally by elements $\{dx,d\log e_i\}$, for $x$ a
   local section of $\O_Y$ and $e_i$ a local generator of $\L_i$, subject to the usual relations among the $dx$,
   plus
   \begin{subequations}\label{E:dlog}
      \begin{equation}
         s_i(e_i)d\log e_i=ds_i(e_i)
      \end{equation}
      \begin{equation}
         \sum d\log e_i=0\text{ when }1\in\O_Y\text{ is written }\otimes e_i
      \end{equation}
   \end{subequations}
   This construction and a subsequent axiomatic development was proposed by Deligne in~\cite{pD88}, and
   independently by Faltings in~\cite{gF90}. In order to deal with the general case ($Y$ no longer assumed to
   have simple normal crossings), it is more convenient to consider, instead of the pairs $(\L_i,s_i)$, the sheaf
   of monoids $\M$ on $Y$ (multiplicatively) generated by $\O^*_Y$ and the $e_i$ together with the map
   $\M\to\O_Y$ given by $\alpha(\otimes e_i^{n_i})=\prod s_i(e_i)^{n_i}$ (replace $d\log e_i$ by $d \log e$ for
   $e\in\M$, and the relations \eqref{E:dlog} by $d\log(ef)=d\log e+d\log f$, $\alpha(e)d\log e=d\alpha(e)$,
   $d\log u=u^{-1}du$ if $u$ is a unit). This is the origin of the notion of logarithmic structure, proposed by
   Fontaine and the speaker, which gave rise to the whole theory beautifully developed by Kato (and Hyodo)
   in~\cite{HK94,kK89a,kK89,kK94b,kK94}.
\end{quotation}

Log structures are natural generalizations to arbitrary schemes of reduced normal crossings divisors on regular
schemes. Here $\M$ is the sheaf of regular functions on the scheme that are invertible away from the divisor,
which is a sheaf of monoids under multiplication. Log structures give rise to differentials with log poles,
crystals and crystalline cohomology with log poles, and similar structures.

In~\cite{kK94}, Kato defines a notion of log regularity and proves that fine saturated log regular schemes behave very
much like toric varieties. If $X$ is such a scheme, then $X$ is normal and Cohen-Macaulay. In addition, Kato develops a
theory of fans and subdivisions, shows how this theory can be used to resolve toric singularities, and identifies the
dualizing complexes of such schemes.

In Chapter~\ref{Ch:regularity} we extend this theory by relaxing the requirement that the monoids be saturated, thereby
relaxing the requirement that the schemes be normal. This is accomplished by casting the theory in terms of log
structures alone. We do not appeal to the traditional theory of cones to keep track of the combinatorics since such an
appeal forces the restriction to saturated monoids. Most of our results are local; that is, they are concerned with the
singularities that occur on these schemes. We define the notion of a toric log regular scheme. Toric varieties are
examples of such schemes, and toric log regular schemes behave like toric varieties in many ways, though they need not
be Cohen-Macaulay. Although most of Kato's methods go through with only minor modifications, there is one significant
change. His induction proof of \cite[(7.2)~Proposition]{kK94} involves a reduction to the case where a chart comes from
a one-dimensional affine semigroup. Once he has reduced to this case, he uses the fact that $\N$ is the only such
moniod. In Lemma~\ref{L:tech}, we introduce a faithful flat descent argument to handle the one-dimensional case.

In Chapter~\ref{Ch:prelim} we collect the algebraic preliminaries will we use in later chapters and the aspects of
basic log geometry that we will emphasize. In Chapter~\ref{Ch:t-flat} we define and study t-flatness.

The notion of t-flatness is due to O.~Gabber. As such, this presentation owes much to him.

\section{Preliminaries}\label{Ch:prelim}

See Kato~\cite{kK89} for an introduction to log schemes. There Kato defines pre-log structures and log structures on
the \'etale site of $X$. However, we will use the Zariski topology throughout this paper. See Niziol~\cite{wN98} for a
brief comparison of log structures on the Zariski and \'etale sites. In Kato~\cite{kK89a}, log structures on locally
ringed spaces are defined and (for the most part) the Zariski topology is used.

\subsection{Fine Log Scheme Basics}

\begin{defn}\cite[Definition~1.2.3]{kK89a}
   Given a pre-log structure $\beta:\M\to\O_X$ on a scheme $X$, the \emph{log structure associated to this pre-log
   structure} $\M^{\alpha}$ is defined to be the colimit of the diagram
   \[
      \xymatrix{
         \beta^{-1}(\O_X^*) \ar[d] \ar[r] & \M \\
         \O_X^* &    }
   \]
   in the category of sheaves of monoids on $X$, equipped with the homomorphism of sheaves of monoids induced by
   $\beta$ and the inclusion $\O_X^*\subseteq\O_X$.
\end{defn}

\begin{rmk}\cite[Section 1.3]{kK89}\label{R:units}
   If $G\xleftarrow{\beta}P\xrightarrow{\gamma}Q$ is a diagram of monoids and $G$ is a group, then the
   colimit of this diagram is $(G\oplus Q)/\sim$, where $\sim$ is the congruence given by $(g,q)\sim(g',q')$ if there
   exist $p,p'\in P$ such that $g+\beta(p)=g'+\beta(p')$ and $q+\gamma(p)=q'+\gamma(p')$. That is, $(g,q)$ and
   $(g',q')$ differ by an element of $P^{gp}$. In particular, we may replace $P$ with any submonoid of
   $(\gamma^{gp})^{-1}(Q)$ that generates $P^{gp}$ without changing the colimit.
\end{rmk}

\begin{xmpl}
   If $P$ is a monoid, $A$ is a ring, and $\beta:P\to A$ is a monoid homomorphism with respect to multiplication on
   $A$, then $\beta$ induces a homomorphism of sheaves of monoids on $\Spec A$ from $P_{\Spec A}$ the constant sheaf on
   $\Spec A$ with stalk $P$ to $\O_{\Spec A}$. We will denote the associated log scheme by $\Spec(P
   \xrightarrow{\beta}A)$.
\end{xmpl}

A coherent log structure $\M$ on a scheme $X$ is integral if and only if locally on $X$, $\M$ is isomorphic to the log
structure associated to the pre-log structure $P_X\to\O_X$ for some finitely generated (cancellative) monoid $P$
(see~\cite[Definition~1.2.6]{kK89a}).

\begin{defn}
   If $\beta:P\to Q$ is a  monoid homomorphism, we say $\beta$ is \emph{local} if $P^*=\beta^{-1}(Q^*)$. Furthermore,
   if $A$ is a ring, $\p\subset A$ is prime, and $\beta:P\to A$ is a monoid homomorphism with repect to multiplication
   on $A$, we say $\beta$ is \emph{local at $\p$} if the composite monoid homomorphism $P
   \xrightarrow{\beta}A\xrightarrow{\text{canonical}}A_{\p}$ is local.
\end{defn}

For the rest of this section, let $P$ be a finitely generated monoid. Let $A$ be a Noetherian ring and let $\beta:P\to
A$ be a monoid homomorphism with respect to multiplication on $A$. Then $\Spec(P\xrightarrow{\beta}A)$ is a Noetherian
fine log scheme. (In fact, every locally Noetherian fine log scheme is locally isomorphic to such a log scheme.) By
Remark~\ref{R:units}, $\Spec(P\xrightarrow{\beta}A)$ is isomorphic to $\Spec(Q\xrightarrow{\beta^{gp}|_Q}A)$ for any
$Q\subseteq P+\beta^{-1}(A^*)$ containing $P$. After we prove the following proposition, we may assume $\beta$ is local
at any particular prime $\p\subset A$ as well.

\begin{prop}\label{P:local}
   If $\p$ is a prime ideal of $A$, then there exists an element $f\in A\setminus\p$ and a finitely generated monoid
   $Q\subseteq P^{gp}$ containing $P$ such that the map $Q\to A_f$ induced by $\beta$ is local at $\p$ and generates
   the same log structure on $\Spec A_f$ as $\beta$.
\end{prop}

\begin{proof}
   Compose $\beta$ with the canonical map $A\to A_{\p}$ to get a map $\widetilde{\beta}:P\to A_{\p}$. Let $X$ be a
   finite generating set for $P$, let $X_0=\{x\in X\mid\widetilde{\beta}(x)\in A_{\p}^*\}$, let $x_0$ be the sum
   of the elements of $X_0$, let $Q$ be the submonoid of $P^{gp}$ generated by $X\cup\{-x_0\}$, and let
   $f=\beta(x_0)$. The rest is straightforward.
\end{proof}

\subsection{Affine Semigroups}

A theorem of Grillet~\cite[Theorem~3.11]{GSR99b} says, ``Let $P$ be a finitely generated monoid. The monoid $P$ is
cancellative, reduced and torsion-free if and only if it is isomorphic to a submonoid of $\N^k$ for some positive
integer $k$.'' We will strengthen this theorem by showing that such a monoid can be embedded in $\N^d$ with $d=\rank
P^{gp}$ in such a way that a complete flag in $P$ is taken to the standard flag in $\N^d$.

\begin{defn}
   A monoid $P$ is said to be \emph{sharp} (or \emph{reduced}) if its unit group is trivial.

   A monoid $P$ is said to be \emph{torsion-free} if its difference group $P^{gp}$ is torsion-free.

   A monoid $P$ is said to be an \emph{affine semigroup} if it is isomorphic to a finitely generated submonoid of
   $\N^k$ for some $k\geq 0$.

   A monoid $S$ is said to be \emph{saturated} (or \emph{normal}) if for any positive integer $n$ and any $p\in
   P^{gp}$, $np\in P$ implies $p\in P$. We call the smallest saturated submonoid of $P^{gp}$ containing $P$
   the \emph{saturation} of $P$ and denote it by $P^{sat}$.
\end{defn}

\begin{prop}\label{P:splitting}
   Let $P$ be a monoid.
   \begin{itemize}
      \item[(1)] If $P$ is sharp and saturated, then $P$ is torsion-free.
      \item[(2)] $\overline{P}$ is sharp.
      \item[(3)] If $P$ is saturated, then $\overline{P}$ is also saturated.
      \item[(4)] If $P$ is saturated, then every localization of $P$ is also saturated.
      \item[(5)] If $P$ is finitely generated and saturated, then $P\cong P^*\bigoplus\overline{P}$.
   \end{itemize}
\end{prop}

\begin{proof}
   (1) We will prove this by contradiction. Suppose $p$ were a nontrivial torsion element of $P^{gp}$ and suppose
   $p$ had order $n$. Then we would have $np=0$. So, $p$ would be in $P$. But, $\N p$ would be a group. So, $p$
   would be a unit in $P$. This contradicts the fact that $P$ is sharp.

   (2) Suppose $\overline{p}\in\overline{P}^*$ and let $p\in P$ be mapped to $\overline{p}$ by the canonical map.
   Since $p$ maps to a unit, there is an element $q\in P$ such that $p+q\in P^*$. Therefore, $q+(-(p+q))$ is the
   inverse of $p$ and $p\in P^*$. Hence, $\overline{p}=0$.

   (3) Suppose $\overline{p}\in\overline{P}^{gp}$ and $n\overline{p}\in\overline{P}$ with $n\in\N$ positive. Let
   $p\in P^{gp}$ be a pre-image of $\overline{p}$. Since $n\overline{p}\in\overline{P}$, $np$ differs from an
   element of $P$ by an element of $P^*$. That is, $np\in P$. Since $P$ is saturated, $p\in P$. Thus,
   $\overline{p}\in\overline{P}$.

   (4) Let $Q=P-S$ with $S$ a submonoid of $P$, let $n$ be a positive integer, and let $q$ be an element of
   $P^{gp}=Q^{gp}$ such that $nq\in Q$. Write $nq=p-s$ with $p\in P$ and $s\in S$. We have $n(q+s)=p+(n-1)s$. So,
   $n(q+s)\in P$. Since $P$ is saturated, $q+s\in P$. So, $q\in Q$.

   (5) By (2), $\overline{P}$ is sharp. By (3), $\overline{P}$ is saturated since $P$ is saturated. By (1),
   $\overline{P}$ is torsion-free. Furthermore, $\overline{P}^{gp}$ finitely generated and hence free. Therefore,
   \[
      \xymatrix@C=0.5cm{
        0 \ar[r] & P^*  \ar[rr] && P^{gp} \ar[rr] && \overline{P}^{gp} \ar[r] & 0 }
   \]
   is split exact. If $p\in P^{gp}$ is mapped into $\overline{P}$ by the righthand map, then $p$ differs from an
   element of $P$ by an element of $P^*$. Hence, $p\in P$ and the pre-image of $\overline{P}$ by the righthand map is
   $P$. That is, any section of $P^{gp}\to\overline{P}^{gp}$ maps $\overline{P}$ into $P$ and $P\cong
   P^*\bigoplus\overline{P}$.
\end{proof}

\begin{rmk}\label{R:torsion-free}
   In Proposition~\ref{P:splitting}, neither (3) nor (4) need be true if the word ``saturated'' is replaced with
   ``torsion-free''. Consider $$P=\langle(1,0),(1,1),(0,2),(0,-2)\rangle\subseteq\Z^2,$$ this monoid is torsion-free but is not saturated. However, $\overline{P}\cong\langle a,b\mid2a=2b\rangle$ is not torsion-free
   (since $a-b\notin\overline{P}$ and $2(a-b)=0$) and $P\ncong P^*\bigoplus\overline{P}$ (since $\overline{P}$ is not
   torsion-free).
\end{rmk}

Consider $V=\Q\otimes_{\Z}P^{gp}$ the $\Q$-vector space generated by the monoid $P$. Let $C(X)$ be the cone over the
subset $X\subseteq P$ in $V$, that is $C(X)=\{\sum_{i=1}^nq_ip_i\in V\mid\forall i,q_i\in\Q_{\geq0},\,p_i\in P\}$. If
$P$ is torsion-free, $P\to V$ is injective so we may freely identify $P$ with its image in $V$. In this case,
$P^{sat}=C(P)\cap P^{gp}$. In particular, if $P$ is torsion-free, $C(P)=C(P^{sat})$.

\begin{defn}
   Let $P$ be a monoid. A submonoid $F\subseteq P$ is said to be a \emph{face} of $P$ if $p+p'\in F$ implies
   $p\in F$.
\end{defn}

Notice that $F\subseteq P$ is a face if and only if $P\setminus F$ is a prime ideal. (We consider $P$ to be a face and
$\emptyset$ to be a prime ideal.)

\begin{prop}\label{P:faces}
   Let $P$ be a finitely generated, torsion-free monoid.
   \begin{itemize}
      \item[(1)] If $F$ is a face of $P$, then $C(F)$ is a face of $C(P)$.
      \item[(2)] If $F$ is a face of $C(P)$, then $P\cap F$ is a face of $P$.
      \item[(3)] If $F$ is a face of $P$, then $P\cap C(F)=F$.
      \item[(4)] If $F$ is a face of $C(P)$, then $C(P\cap F)=F$.
   \end{itemize}
   This establishes a bijective correspondence between the faces of $P$ and the faces of $C(P)$.
\end{prop}

\begin{proof}
   Evidently, whenever $X\subseteq P$ and $x\in C(X)$, there exists a positive integer $n$ such that $nx$ is in the
   monoid generated by $X$.

   (1) Suppose $x$ and $y$ are elements of $C(P)$ such that $x+y\in C(F)$. Let $n_1$ be a positive integer such that
   $n_1x\in P$, let $n_2$ be a positive integer such that $n_2y\in P$, let $n_3$ be a positive integer such that
   $n_3(x+y)\in F$, and let $m$ be the least common multiple of $n_1$, $n_2$ and $n_3$. We have $mx\in P$, $my\in P$
   and $mx+my=m(x+y)\in F$. Since $F$ is a face of $P$, $mx\in F$ and $my\in F$. So, $x\in C(F)$ and $y\in C(F)$.

   (2) Suppose $p$ and $p'$ are elements of $P$ such that $p+p'\in P\cap F$. Since $F$ is a face, both $p$ and $p'$ are
   in $F$. So, $p\in P\cap F$ and $p'\in P\cap F$.

   (3) Suppose $p$ is an element of $P\cap C(F)$. Let $n$ be a positive integer such that $np\in F$. Since $(n-1)p+p\in
   F$ and $F$ is a face, $p\in F$. Evidently, $F\subseteq P\cap C(F)$.

   (4) If $x\in C(P\cap F)$, then for some positive integer $n$, $nx\in P\cap F$. In particular, $nx\in F$. Since $F$
   is a face, $x\in F$. If $x\in F$, then for some positive integer $n$, $nx\in P\cap F$. In particular, $nx\in C(P\cap
   F)$. Since $C(P\cap F)$ is a face, $x\in C(P\cap F)$.
\end{proof}

In light of this, when $P$ is finitely generated and torsion-free, we freely speak of edges (resp.\ facets etc.)
meaning faces whose corresponding faces in $C(P)$ are edges (resp.\ facets etc.).

\begin{defn}
   A sequence of faces $F_0\subsetneq F_1\subsetneq F_2\subsetneq\cdots\subsetneq F_r$ in a monoid $P$ is said to
   be a \emph{flag} in $P$. A flag in $P$ is said to be \emph{complete} if $F_0=P^*$, $F_r=P$ and whenever
   $F$ is a face of $P$ lying between $F_{i-1}$ and $F_i$, $F=F_{i-1}$ or $F=F_i$.
\end{defn}

The complements in $P$ of the faces in a flag form a chain of prime ideals. If the flag is complete, the chain of
primes is saturated; that is, there is no prime properly between consecutive primes of the chain. Furthermore, if
$\varphi:P\to Q$ is a homomorphism of monoids and $G$ is a face of $Q$, $\varphi^{-1}(G)$ is a face of $P$. We will say
$\varphi$ takes the flag $G_0\subsetneq G_1\subsetneq G_2\subsetneq\cdots\subsetneq G_r$ in $Q$ to the flag
$F_0\subsetneq F_1\subsetneq F_2\subsetneq\cdots\subsetneq F_r$ in $P$ via pullback if $\varphi^{-1}(G_i) = F_i$ for
all $0\leq i\leq r$.

Let $\mathbf{e}_i\in\N^d$ be the element with a $1$ in the $i$th position and zeroes elsewhere. We call
$\{0\}\subset\langle\mathbf{e}_1\rangle\subset\langle\mathbf{e}_1,\mathbf{e}_2\rangle\subset\cdots\subset\N^d$ the
standard flag in $\N^d$.

\begin{thm}\label{T:inclusion}
   Let $P$ be a sharp finitely generated torsion-free monoid and let $\{0\}=F_0\subsetneq F_1\subsetneq
   F_2\subsetneq\cdots\subsetneq F_d$ be a complete flag in $P$. Then there exists an inclusion $\varphi:P\to\N^d$
   taking the standard flag in $\N^d$ to the given flag in $P$ and inducing an isomorphism
   $\varphi^{gp}:P^{gp}\to\Z^d$.
\end{thm}

\begin{proof}
   Without loss of generality we may assume $P$ is saturated, since the faces of $P$ correspond bijectively to the
   faces of $C(P)=C(P^{sat})$ according to Proposition~\ref{P:faces}. We proceed by induction on $d$. When $d=0$, the
   theorem is trivial.

   Suppose we have a map $\widetilde{\varphi}:F_{d-1}\to\N^{d-1}$ that sends the flag $\{0\}=F_0\subsetneq
   F_1\subsetneq F_2\subsetneq\cdots\subsetneq F_{d-1}$ to the standard flag in $\N^{d-1}$ and inducing an isomorphism
   $\widetilde{\varphi}^{gp}:F_{d-1}^{gp}\to\Z^{d-1}$. Let $\sigma$ be a splitting of the inclusion $F_{d-1}^{gp}\to
   P^{gp}$. Such a splitting $\sigma$ exists since $F_{d-1}$ is a face and $P$ is saturated by
   Proposition~\ref{P:splitting}. Now let $\psi_1:P\to\Z^{d-1}$ be the map
   \[
      \xymatrix{
         P \ar[r] & P^{gp} \ar[r]^{\sigma} & F_{d-1}^{gp} \ar[r] & \Z^{d-1}.   }
   \]
   That is, inclusion into $P^{gp}$ followed by $\sigma$ followed by the isomorphism $\widetilde{\varphi}^{gp}$.
   Notice that if $p\in F_{d-1}$, then the coordinates of $\psi_1(p)$ are non-negative.

   Let $\p$ be the complement of $F_{d-1}$. Since $F_{d-1}$ is a facet of $P$, $\p$ is a height one prime of $P$.
   $\overline{P_{\p}}\cong\N$ since $P$ was assumed to be saturated. Let $\psi_2:P\to\N$ be be the map
   \[
      \xymatrix{
         P \ar[r] & P_{\p} \ar[r] & \overline{P_{\p}} \ar[r]^{\cong} & \N.   }
   \]
   Notice that $\psi_2(p)=0$ if and only if $p\in F_{d-1}$. Now I claim $\psi=(\psi_1, \psi_2):P\to\Z^{d-1}\times\N$ is
   injective. If $\psi(p)=\psi(p')$, then their last coordinates are equal and $p-p'\in F_{d-1}^{gp}$. Since the map
   $F_{d-1}^{gp}\to\Z^{d-1}$ above is an isomorphism, $\psi$ is injective.  Furthermore, since some element of $P$ maps
   to an element of $\Z^d$ whose last coordinate is one and the map $F_{d-1}^{gp}\to\Z^{d-1}$ above is an isomorphism,
   $\psi^{gp}$ is an isomorphism.

   Recall that the unique minimal generating set of an affine semigroup is called its Hilbert basis. For every element
   $p$ of the Hilbert basis of $P$, let $$(p_1,p_2,\ldots,p_d)=\psi(p).$$ For $1\leq i\leq d$, choose $n_i\in \N$ such
   that $p_i+n_ip_d\geq0$ for every Hilbert basis element $p\in P$, and let $\theta$ be the automorphism of $\Z^d$
   given by
   \[
      \left[
         \begin{matrix}
            1      & 0      & 0      & \cdots & 0      & n_1     \\
            0      & 1      & 0      & \cdots & 0      & n_2     \\
            0      & 0      & 1      & \cdots & 0      & n_3     \\
            \vdots & \vdots & \vdots & \ddots & \vdots & \vdots  \\
            0      & 0      & 0      & \cdots & 1      & n_{d-1} \\
            0      & 0      & 0      & \cdots & 0      & 1
         \end{matrix}
      \right]
   \]
   Take $\varphi=\theta\circ\psi$.
\end{proof}

\begin{xmpl}
   Let $P=\langle(0,2),(1,0),(2,-2)\rangle\subset\N^2$ and consider the flag
   $\{0\}\subset\langle(0,2)\rangle\subset P$. Suppose we already have
   $\widetilde{\varphi}:\langle(0,2)\rangle\to\N$, given by $(0,2)\mapsto1$. Let $\psi_1:P\to\Z$ be given by
   $(0,2)\mapsto1$, $(1,0)\mapsto 0$ and $(2,-2)\mapsto-1$. $\psi_2:P\to\N$ is given by $(0,2)\mapsto0$,
   $(1,0)\mapsto1$ and $(2,-2)\mapsto2$. Therefore, $\psi:P\to\Z\times\N$ is given by $(0,2)\mapsto(1,0)$,
   $(1,0)\mapsto(0,1)$ and $(2,-2)\mapsto(-1,2)$. Let $n_1=1$, then
   \[
      \theta =
      \left[
         \begin{matrix}
            1 & 1 \\
            0 & 1
         \end{matrix}
      \right]
   \]
   So, $\varphi$ is given by $(0,2)\mapsto(1,0)$, $(1,0)\mapsto(1,1)$ and $(2,-2)\mapsto(1,2)$.
\end{xmpl}

\subsection{Toriodal Log Schemes}

\begin{defn}
   We say a log structure $\M$ on $X$ is \emph{torsion-free} if $\overline{\M}$ is a sheaf of torsion-free monoids.
   We say a log structure $\M$ on $X$ is \emph{toriodal} if it is fine and torsion-free. We say a log scheme
   $(X,\M)$ is \emph{toriodal} if $\M$ is toroidal.
\end{defn}

\begin{ass}
   If $P$ is a sharp finitely generated torsion-free monoid, $A$ is a Noetherian ring, and $\beta:P\to A$ is amonoid
   homomorphism with respect to multiplication on $A$, $\Spec(P\xrightarrow{\beta}A)$ need not be toroidal. Let $k$ be
   a field and consider
   \[
      \beta:\langle(1,0),(1,1),(0,2)\rangle\to k[x,y,z]/(x^2z-y^2)
   \]
   given by $(1,0)\mapsto x$, $(1,1)\mapsto y$, and $(0,2)\mapsto z$. Here $\overline{M}_{(x,y)}\cong\langle
   a,b\mid2a=2b\rangle$. See Remark~\ref{R:torsion-free}.
\end{ass}

\begin{prop}\label{P:canonical_chart}
   If $(X,\M)=\Spec(P\xrightarrow{\beta}A)$ is a fine log scheme and $\beta:P\to A$ is local at
   $x=\p$, then $\overline{\M}_x\cong\overline{P}$. Furthermore, if
   \[
      \xymatrix@C=0.5cm{
         0 \ar[r] & P^* \ar[rr] && P^{gp} \ar[rr] && \overline{P}^{gp} \ar[r] & 0 }
   \]
   is split exact, then there is a monoid
   homomorphism $\sigma:\overline{P}\to P$ such that the morphism
   $\Spec(\overline{P}\xrightarrow{\sigma\circ\beta}A)\to\Spec(P\xrightarrow{\beta}A)$ induced by $\sigma$ is an
   isomorphism.
\end{prop}

\begin{proof}
   Since the question of whether $\overline{\M}_x$ is isomorphic to $\overline{P}$ is local, we may assume
   $x$ is contained in every irreducible component of $X$. In particular, we may assume any neighborhood $U$ of
   $x$ is connected and $\Gamma(U,P_X)=P$. Let $\overline{\beta}:P\to A_{\p}$ be the composition of $\beta$ and
   the canonical map $A\to A_{\p}$. By Remark~\ref{R:units}, $\M_x\cong A_{\p}^*\oplus
   P/\{\overline{\beta}(u)^{-1}\oplus u\mid u\in P^*\}$. We may now take the quotient $\M_x/\M_x^*$ to
   form $\overline{\M}_x$. Evidently, $\overline{\M}_x\cong(A_{\p}^*\oplus P)/(A_{\p}^*\oplus P^*)$. But,
   $\overline{P}\cong(A_{\p}^*\oplus P)/(A_{\p}^*\oplus P^*)$. The rest is straightforward. Any splitting of the
   exact sequence restricts to such a $\sigma:\overline{P}\to P$.
\end{proof}

\begin{cor}\label{C:chart}
   Let $(X,\M)$ be a fine log scheme and let $x$ be a point on $X$ such that
   $\overline{\M}_x$ is a torsion-free monoid. Then, there exists an open neighborhood $U$ of $x$ and a
   homomorphism $\beta:\overline{\M}_x\to\O_X(U)$ such that $\beta$ induces $\alpha|_U$ and the composition of
   the induced map $\overline{\M}_x\to\M_X(U)$ with restriction to the stalk $\M_x$ and the canonical map
   $\M_x\to\overline{\M}_x$ is the identity.
\end{cor}

\subsection{The Prime Filtration Theorem}\label{S:prime_filtration}

\begin{thm}\emph{[Prime Filtration Theorem]}\label{T:prime_filtration}
   Let $P$ be a finitely generated torsion-free monoid, let $A=\bigoplus_{p\in P}A_p$ be a Noetherian $P$-graded
   ring, and let $E=\bigoplus_{p\in P^{gp}}E_p$ be a finitely generated $P^{gp}$-graded $A$-module. Then, $E$ has
   a filtration
   \[
      0=E_0\subset E_1\subset\cdots\subset E_n=E
   \]
   with each $E_i/E_{i-1}\cong A/\p_i$ for some homogeneous prime ideal $\p_i\subset A$.
\end{thm}

\begin{proof}
   By the usual prime filtration theorem \cite[Proposition~3.7]{dE95}, $E$ has a filtration
   \[
      0=E_0\subset E_1\subset\cdots\subset E_n=E
   \]
   with each $E_i/E_{i-1}\cong A/\p_i$ for some prime ideal $\p_i$. By induction on $n$, it suffices to prove
   that $\p_1$ is homogeneous.

   Since $A/\p_1\subseteq E$, $\p_1$ is an associated prime of $E$. Hence, there is a least positive integer $s$
   such that there exists an element $e$ of $E$ with annihilator $\p$ and $P^{gp}$-homogeneous elements $e_i$ of
   $E$ such that $e=\sum_{i=1}^s e_i$. By the results of \S~3 of Gilmer~\cite{rG84}, there is a total order
   $\prec$ compatible with the group structure on $P^{gp}$ since $P^{gp}$ is torsion-free. For each $i$, let
   $q_i$ be the degree of $e_i$. We may assume $q_1\prec q_2\prec\ldots\prec q_s$. Let $f$ be an element of
   $\p_1$. Write $f=\sum_{i=1}^mf_i$ with each $f_i$ $P$-homogeneous of degree $p_i$ and $p_1\prec
   p_2\prec\ldots\prec p_m$. We want to show that each homogeneous component $f_i$ of $f$ is in $\p_1$. By
   induction on $m$, it suffices to prove $f_1$ is in $\p_1$.

   Notice that
   \[
      0=fe=f_1e_1+(\text{homogeneous terms of greater degree}).
   \]
   So, $f_1e_1$ is zero and we are done if $s=1$, that is, if $e$ is $e_1$. Notice that $f_1e$ is
   $\sum_{i=2}^sf_1e_i$ and has fewer nonzero homogeneous terms than $e$. Since $s$ was minimal, the annihilator
   of $f_1e$ properly contains $\p_1$. Let $f'$ be an element of the annihilator of $f_1e$ that is not contained
   in $\p_1$. We know $f'f_1$ is in $\p_1$, $f'$ is not in $\p_1$, and $\p_1$ is prime. Therefore, $f_1$ is in
   $\p_1$.
\end{proof}

\begin{rmk}
   We cannot delete the torsion-free assumption in the above theorem: Let $R$ be any ring, let $a$ and $b$ be
   distinct elements of $P$ such that $na=nb$ for some positive integer $n$, and let $n_0$ be the least such
   positive integer. Now $t^a-t^b$ is an element of $R[P]$ that is annihilated by
   $\sum_{m=0}^{n_0-1}t^{(n_0-m-1)a+mb}$. But, $t^a-t^b$ cannot be annihilated by any nonzero homogeneous
   element of $R[P]$ because $P$ is cancellative.
\end{rmk}

\subsection{Combinatorial $\Z[P]$-modules}\label{S:combinatorial}

\begin{lemma}\label{L:homogeneous_prime}
   Let $P$ be a finitely generated torsion-free monoid. If $\p$ is a $P$-homogeneous prime of $\Z[P]$, then
   $\p=(q)+(K)=(q,K)$ where $(q)\subset\Z$ and $K\subset P$ are prime.
\end{lemma}

\begin{proof}
   Fix an arbitrary $P$-homogeneous prime $\p$ of $\Z[P]$. The intersection $\p\cap\Z$ must be prime. Let
   $(q)=\p\cap\Z$. If $\p=q\Z[P]$ we are done, take $K$ to be the empty ideal. If $\p\neq q\Z[P]$, pick an
   arbitrary $P$-homogeneous element $nt^p$ of $\p\setminus q\Z[P]$ with $n\in\Z$ and $p\in P$. In particular,
   $n\notin (q)=\p\cap\Z$. So, $n\notin\p$. Hence, $t^p\in\p$ since $\p$ is prime. Let $K=\{p\in P\mid
   t^p\in\p\}$ be the set of all such $p$. We have established $\p$ is generated by $q$ and $(K)$. Now, it
   suffices to prove $K$ is a prime ideal of $P$. If $k$ is in $K$ and $p$ is in $P$, then $t^{p+k}=t^p\cdot t^k$
   is in $\p$ since $\p$ is an ideal. So, $K$ is an ideal. If $K$ were not prime, there would exist $p$ and $p'$
   in $P$ such that $p+p'$ would be in $K$ while neither $p$ nor $p'$ was in $K$. That is, there would exist
   $t^p$ and $t^{p'}$ in $\Z[P]$ such that $t^p\cdot t^{p'}$ would be in $\p$ while neither $t^p$ nor $t^q$ was
   in $\p$ and $\p$ would not be prime. Hence, $K$ is prime.
\end{proof}

\begin{defn}
   Let $P$ be a monoid. We say $A$ is a \emph{sub-$P$-set} of $P^{gp}$ if $A$ is a subset of $P^{gp}$ and $A$ is
   closed under the action of $P$ on $P^{gp}$ given by addition. That is, $P+A=A$. We say a sub-$P$-set of
   $P^{gp}$ is a \emph{fractional ideal} if there exists an element $p$ of $P$ such that $p+A$ is contained in
   $P$. A \emph{combinatorial $\Z[P]$-module} is one isomorphic to a quotient of $\Z[P]$-submodules $(A)/(B)$ of
   $\Z[P^{gp}]$, for some sub-$P$-sets $A$ and $B$ of $P^{gp}$. If $K$ is an ideal of $P$ and our module is also
   annihilated by $(K)$, we will say it is a \emph{combinatorial $\Z[P]/(K)$-module}.
\end{defn}

\begin{cor}\label{C:filtration}
   If $P$ is a finitely generated torsion-free monoid and $N$ is a finitely generated combinatorial
   $\Z[P]$-module, then $N$ has a filtration
   \[
      0=N_0\subset N_1\subset\cdots\subset N_n=N
   \]
   with each $N_i/N_{i-1}\cong\Z[P]/(K)$ for some prime ideal $K\subset P$.
\end{cor}

\begin{proof}
   We may assume $N=(A)/(B)$ where $A$ and $B$ are fractional ideals of $P$ since $P$ is a finitely generated
   torsion-free monoid. By Theorem~\ref{T:prime_filtration}, we know $N$ has a filtration
   \[
      0=N_0\subset N_1\subset\cdots\subset N_n=N
   \]
   with each $N_i/N_{i-1}\cong\Z[P]/\p_i$ for some $P$-homogeneous prime ideal $\p_i$. In particular, for
   $\p_1$ we may take any associated prime of $N$. By induction on $n$, it suffices to prove that $N$ has
   an associated prime of the form $(K)$ for some prime $K\subset P$.

   By Lemma~\ref{L:homogeneous_prime}, every associated prime of a finitely generated combinatorial
   $\Z[P]$-module is of the form $(q,K)$ for some primes $K\subset P$ and $(q)\subset\Z$. Evidently, each
   combinatorial $\Z[P]$-module is a free abelian group and on any combinatorial $\Z[P]$-module every integer acts
   injectively. So, our prime is of the form $(K)$ for some prime $K\subset P$.
\end{proof}

\begin{prop}\label{P:inclusion}
   If $P\to Q$ is an inclusion of finitely generated torsion-free monoids, $K$ is an ideal of $Q$, $K'=K\cap P$, and
   $N$ is a combinatorial $\Z[Q]/(K)$-module, then $N$ is a direct sum of combinatorial $\Z[P]/(K')$-modules.
\end{prop}

\begin{proof}
   Write $N=(A)/(B)$ with $A$ and $B$ sub-$Q$-sets of $Q^{gp}$ such that $A+K\subseteq B$. Pick a set $C$ of coset
   representatives for $Q^{gp}/P^{gp}$. For each $c\in C$, let $A_c=\{q-c\mid q\in A\text{ and }q-c\in P^{gp}\}$ and
   let $B_c=\{q-c\mid q\in B\text{ and }q-c\in P^{gp}\}$. For each $c\in C$, $A_c$ and $B_c$ are sub-$P$-sets of
   $P^{gp}$ such that $A_c+K'\subseteq B_c$. Therefore, $(A_c)/(B_c)$ is a combinatorial $\Z[P]/(K')$-module.
   Furthermore, $(A)/(B)\cong\bigoplus_{c\in C}(A_c)/(B_c)$ as $\Z[P]/(K')$-modules.
\end{proof}

\subsection{Lorenzon's Algebra}

Let $(X,\M)$ be a fine log scheme. For every local section $\overline{m}$ of $\overline{\M}$, the
pre-image of $\overline{m}$ along the canonical map $\M\to\overline{\M}$ is an $\O_X^*$-torsor. So, the canonical
map $\M\to\overline{\M}$ is an $\overline{\M}$-indexed family of $\O_X^*$-torsors. To each $\O_X^*$-torsor $\E$,
we associate the contracted product $\E\wedge_{\O_X^*}\O_X$, where $\E\wedge_{\O_X^*}\O_X$ is the quotient of the
product $\E\times\O_X$ by the equivalence relation $\sim$, where $(eu,f)\sim(e,uf)$ whenever $e$, $u$ and $f$ are
respectively local sections of $\E$, $\O_X^*$ and $\O_X$. Each $\E\wedge_{\O_X^*}\O_X$ is an invertible sheaf. The
association of $\E\wedge_{\O_X^*}\O_X$ to $\E$ and the $\overline{\M}$-indexed family of $\O_X^*$-torsors
$\M\to\overline{\M}$ together yield an $\overline{\M}$-indexed family of invertible sheaves.

Furthermore, if $\overline{m}$ and $\overline{m}'$ are local sections $\overline{\M}$ with corresponding invertible
sheaves $\L_{\overline{m}}$ and $\L_{\overline{m}'}$, then $\L_{\overline{m}+\overline{m}'}$ is isomorphic to
$\L_{\overline{m}}\mathbf{\otimes}\L_{\overline{m}'}$ (see Lorenzon~\cite{pL01}). This family of invertible sheaves is
a ring object in the comma category of sheaves of sets on $X$ over $\overline{\M}$, $\Sh(X)/\overline{\M}$, and the
part indexed by the identity section of $\overline{\M}$ is $\O_X$. Lorenzon calls this $\overline{\M}$-indexed
$\O_X$-algebra the canonical algebra of the log scheme. Let the algebra $\A_X$ be the direct sum over sections of
$\overline{\M}$ of these invertible sheaves with multiplication given by tensor product over $\O_X$ extended by
$\O_X$-linearity. Lorenzon calls this algebra the algebra induced on $X$ by the canonical algebra of the log scheme.
The canonical algebra consists of the homogeneous pieces of our $\overline{\M}$-graded algebra.

\begin{ass}
   Lorenzon uses $(\A_X)_*$ to denote the algebra we denote by $\A_X$ and he uses $\A_X$ to denote the canonical
   algebra of the log scheme.
\end{ass}

It happens that $\A_X$ is the $\O_X$-algebra $\O_X[\M]/(t^m-\alpha(m))_{m\in\M^*}$. For an ideal sheaf $\K$ of $\M$,
let $(\K)$ be the image in $\A_X$ of the ideal sheaf $\K$ along the canonical homomorphism $\M\to\O_X[\M]\to\A_X$.

\begin{xmpl}($\P^1$ with a marked point)
   Let $k$ be a field. Consider the point $p=(x)$ on $\P^1=\Proj k[x,y]$. Let $\M$ be the log structure that is trivial
   on the open subset $\P^1\setminus\{p\}$ and induced by the $\N\to k[\frac{x}{y}],\,1\mapsto\frac{x}{y}$ on the
   affine open subscheme $U=\Spec k[\frac{x}{y}]$ containing $p$. Here $\overline{\M}$ is the skyscraper sheaf with
   stalk $\N$ at $p$ and the $\O_{\P^1}^*$-torsor associated to the natural number $n$ is generated by
   $(\frac{x}{y})^n$ on $U$ and by $1$ on $\P^1\setminus\{p\}$. So, the invertible sheaf associated to the natural
   number $n$ has local basis $(\frac{x}{y})^n$ on $U$ and $1$ on $\P^1\setminus\{p\}$. That is, the invertible sheaf
   associated to the natural number $n$ is $\L(-p)^{\otimes n}$ and $\A_{\P^1}$ is the $\O_{\P^1}$-algebra
   $\bigoplus_{n\in\N}\L(-p)^{\otimes n}$.
\end{xmpl}

\begin{prop}\label{P:A_chart}
   Let $(X,\M)=\Spec(P\xrightarrow{\beta}A)$ be a fine log scheme and let $\beta:P\to A$ be local
   at $x=\p$, then
   \[
      \A_{X,x}\cong\O_x[P]/(t^p-\beta(p))_{p\in P^*}.
   \]
\end{prop}

\begin{proof}
   To prove this, we write
   $\M_x\cong\O_{X,x}^*\oplus P/\{\beta(u)^{-1}\oplus u\mid u\in P^*\}$ and interchange the order in which quotients
   are taken. Instead of first applying the congruence that forms $\M_x$ from $\O_{X,x}^*\oplus P$, we first identify
   the copy of $\O_{X,x}^*$ in $\O^*_{X,x}\oplus P$ with the copy of $\O_{X,x}^*$ in $\O_{X,x}$:
   \begin{align*}
      \A_{X,x}  &=\O_{X,x}[\M_x]/(t^m-\alpha(m))_{m\in \M_x^*} \\
            &=\O_{X,x}[\O_{X,x}^*\oplus P]/
               (t^{(1,p)}-t^{(\beta(p),0)},t^{(u,0)}-u)_{p\in P^*,u\in \O_{X,x}^*} \\
            &=\O_{X,x}[P]/(t^p-\beta(p))_{p\in P^*}
   \end{align*}
\end{proof}

In particular, $\A_{X,x}\cong\O_{X,x}[\overline{\M}_x]$ when $\overline{\M}_x$ is torsion-free by
Proposition~\ref{P:canonical_chart}.

\section{t-Flatness}\label{Ch:t-flat}

\subsection{First Properties of t-Flatness}

\begin{defn}
   Let $(X,\M)$ be an fine log scheme, let $x$ be a point of $X$, let $\K$ be an ideal of
   $\M$ and let $\F$ be a sheaf of $\O_X$-modules such that $(\K)_x\F_x=0$. We say $\F$ is \emph{$\M$-flat
   relative to $\K$ at $x$} if for all ideals $\J$ of $\M$ containing $\K$ we have
   \[
      \Tor_1^{\A_{X,x}/(\K)_x}(\A_{X,x}/(\J)_x,\F_x)=0
   \]
   We say $\F$ is $\M$-flat at $x$ if $\F$ is $\M$-flat relative to the constant ideal sheaf with empty
   stalks at $x$.
\end{defn}

\begin{defn}
   Let $P$ be a finitely generated monoid, let $K\subseteq P$ be an ideal, and let $E$ be a nonzero
   $\Z[P]/(K)$-module. We say $E$ has \emph{t-flat dimension $d$ relative to $K$} if
   \[
      d=\sup\{i\mid\exists J\subseteq P\text{ containing }K,\Tor_i^{\Z[P]/(K)}(\Z[P]/(J),E)\neq0\}.
   \]
   If $E=0$, we say $E$ has t-flat dimension $0$ relative to $K$. We say $E$ is \emph{t-flat relative to $K$} if
   \[
      \Tor_1^{\Z[P]/(K)}(\Z[P]/(J),E)=0
   \]
   for all ideals $J\subseteq P$ containing $K$. If $E$ is t-flat relative to $\emptyset$, we simply say $E$ is
   t-flat.
\end{defn}

Later, Theorem~\ref{T:local_criterion}, we will prove $E$ is t-flat relative to $K$ if and only if $E$ has t-flat
dimension $0$ relative to $K$.

\begin{prop}
   Let $P$ be a finitely generated torsion-free monoid, let $A$ be a Noetherian ring, let $\p\subset A$ be prime, let
   $\beta:P\to A$ be a monoid homomorphism with respect to multiplication on $A$, let
   $(X,\M)=\Spec(P\xrightarrow{\beta}A)$, and let $x$ be the point on $X$ corresponding to $\p$.
   Suppose $\overline{\M}_x$ is torsion-free and let $\F$ be an $\O_X$-module. We consider $\F_x$ to be a
   $\Z[P]$-module along the map $\Z[P]\to\O_{X,x}$ induced by $\beta$. If $K\subseteq P$ is an ideal such that $(K)$
   annihilates $\F_x$, then
   \[
      \Tor_i^{\Z[P]/(K)}(\Z[P]/(J),\F_x)\cong\Tor_i^{\A_{X,x}/(K)_x}(\A_{X,x}/(J)_x,\F_x),\,\forall i\geq0
   \]
   for all ideals $J$ of $P$ containing $K$.
\end{prop}

\begin{proof}
   Since both of these modules are $\O_{X,x}$- modules, both sides are zero if $J\nsubseteq\beta^{-1}(\p)$ and we may
   assume $A=\O_{X,x}$. By Proposition~\ref{P:canonical_chart}, we may assume $P\cong\overline{\M}_x$. In
   particular, we may assume $\A_{X,x}=A[P]$ by Proposition~\ref{P:A_chart}. Fix ideals $K\subseteq J\subseteq P$. Let
   \[
      \mathbf{F}_{\centerdot}:\qquad\xymatrix@C=0.5cm{
         \cdots \ar[rr] && F_2 \ar[rr] && F_1 \ar[rr] && F_0 \ar[rr] && \Z[P]/(J) \ar[r] & 0 }
   \]
   be a $P^{gp}$-graded free resolution of the $\Z[P]/(K)$-module $\Z[P]/(J)$. Since $\Z[P]/(J)$ is a free $\Z$-module,
   \begin{multline*}
      \mathbf{F}'_{\centerdot}:\qquad\xymatrix@C=0.5cm{\cdots \ar[rr] && F_2\otimes_{\Z}A \ar[rr] && F_1\otimes_{\Z}A} \\
      \xymatrix@C=0.5cm{\ar[rr] && F_0\otimes_{\Z}A \ar[rr] && A[P]/(J) \ar[r] & 0}
   \end{multline*}
   is a $P^{gp}$-graded free resolution of $A[P]/(J)$ as a $A[P]/(K)$-module.
   Furthermore,
   \begin{align*}
      F_i&\otimes_{\Z[P]/(K)}\F_x \cong F_i\otimes_{\Z[P]/(K)}A[P]/(K)\otimes_{A[P]/(K)} \F_x \\
         &\cong F_i\otimes_{\Z[P]/(K)}\Z[P]/(K)\otimes_{\Z}A \otimes_{A[P]/(K)}\F_x \\
         &\cong F_i\otimes_{\Z}A\otimes_{A[P]/(K)}\F_x
   \end{align*}
   Therefore, $\Tor_i^{\Z[P]/(K)}(\Z[P]/(J),\F_x)$, the $i$th cohomology module of
   \[
      \mathbf{F}_{\centerdot}\otimes_{\Z[P]/(K)}\F_x,
   \]
   and $\Tor_i^{A[P]/(K)}(A[P]/(J),\F_x)$, the $i$th cohomology module of
   \[
      \mathbf{F}'_{\centerdot}\otimes_{A[P]/(K)}\F_x,
   \]
   are isomorphic for all $i\geq0$.
\end{proof}

\begin{cor}\label{C:comparison}
   Let $P$ be a finitely generated torsion-free monoid, let $A$ be a Noetherian ring, let $\p\subset A$ be prime, let
   $\beta:P\to A$ be a monoid homomorphism with respect to multiplication on $A$, let
   $(X,\M)$ be the log scheme $\Spec(P\xrightarrow{\beta}A)$, and let $x$ be the point on $X$ corresponding to $\p$.
   Suppose $\overline{\M}_x$ is torsion-free and let $\F$ be an $\O_X$-module. We consider $\F_x$ to be a
   $\Z[P]$-module along the map $\Z[P]\to\O_{X,x}$ induced by $\beta$. If $K\subseteq P$ is an ideal such that $(K)$
   annihilates $\F_x$ and $\K\subseteq\M$ is the ideal generated by the image of $K$, then $\F$ is $\M$-flat
   relative to $\K$ at $x\in X$ if and only if $\F_x$ is t-flat relative to $K$.
\end{cor}

\subsection{Local Criterion for t-Flatness}

We will model the proof of our Local Criterion for t-Flatness on the proof of the Local Criterion for Flatness found
in Matsumura~\cite{hM86}.

Recall that a $A$-module $E$ is said to be $I$-adically ideal-separated if $\a\otimes E$
is separated for the $I$-adic topology for every finitely generated ideal $\a$ of $A$. In particular, if $B$ is a
Noetherian $A$-algebra, $IB\subseteq\rad(B)$, and $E$ is a finitely generated $B$-module, then $E$ is $I$-adically
ideal-separated.

\begin{thm}\emph{\cite[Theorem~22.3]{hM86}}
   Let $I$ be an ideal of a ring $A$ and let $M$ be a $A$-module. Set $A_n=A/I^{n+1}$ for  each integer $n\geq0$,
   $M_n=M/I^{n+1}M$ for each integer $n\geq0$, $\gr(A)=\bigoplus_{n\geq0}I^n/I^{n+1}$, and
   $\gr(M)=\bigoplus_{n\geq0}I^nM/I^{n+1}M$. Let
   \[
      \mu_n:(I^n/I^{n+1})\otimes_{A_0}M_0\to I^nM/I^{n+1}M
   \]
   be the standard map for each $n\geq0$, and let
   \[
      \mu:\gr(A)\otimes_{A_0}M_0\to\gr(M)
   \]
   be the direct sum of the $\mu_n$.

   In the above notation, suppose that one of the following two conditions is satisfied: ($\alpha$) $I$ is a nilpotent
   ideal or ($\beta$) $A$ is a Noetherian ring and $M$ is $I$-adically ideal-separated. Then the following conditions
   are equivalent,
   \begin{itemize}
      \item[(1)] $M$ is flat over $A$;
      \item[(2)] $\Tor_1^A(N,M)=0$ for every $A_0$-module $N$;
      \item[(3)] $M_0$ is flat over $A_0$ and $I\otimes_AM=IM$;
      \item[(3')] $M_0$ is flat over $A_0$ and $\Tor_1^A(A_0,M)=0$;
      \item[(4)] $M_0$ is flat over $A_0$ and $\mu_n$ is an isomorphism for every $n\geq0$;
      \item[(4')] $M_0$ is flat over $A_0$ and $\mu$ is an isomorphism;
      \item[(5)] $M_n$ is flat over $A_n$ for every $n\geq0$.
   \end{itemize}
   In fact, the implications
   (1)$\Rightarrow$(2)$\Leftrightarrow$(3)$\Leftrightarrow$(3')$\Rightarrow$(4)$\Rightarrow$(5) hold without the
   assumption on $M$.
\end{thm}

Let $P$ be a finitely generated torsion-free monoid. We consider the case where $A=\Z[P]/(K)$, $I=(P^+)/(K)$ and one of
the following two conditions holds: there exists a Noetherian $A$-algebra $B$ with $IB\subseteq\rad(B)$, and $E$ is a
finitely generated $B$-module or $I$ is nilpotent.

Let $d$ be the dimension of $\overline{P^{sat}}$ and let $\varphi:\overline{P^{sat}}\to\N^d$ be an inclusion as in
Proposition~\ref{T:inclusion}. Pick $d$ $\Q$-linearly independent positive real numbers
$$\{\gamma_1,\gamma_2,\ldots,\gamma_d\},$$ let the monoid homomorphism $\psi:\N^d\to\R$ be given by
$$(n_1,n_2,\ldots,n_d)\mapsto\sum_{i=1}^dn_i\gamma_i,$$ let the monoid homomorphism $\nu:P\to\R$ be the composition
\[
   \xymatrix{P \ar[r]^{\text{can.}} & \overline{P^{sat}} \ar[r]^{\varphi} & \N^d \ar[r]^{\psi} & \R},
\]
and let $\Gamma$ be the image of $\nu$. We order $\Gamma$ with the order induced by the standard ordering of $\R$.
Notice that $\Gamma$ is well ordered.

Set $I_{\gamma}=(t^p)_{\gamma\leq\nu(p),\,p\in P\setminus K}$ for each $\gamma\in\Gamma$, $K_{\gamma}^+=\{p\in
P\mid\nu(p)>\gamma\}$ for each $\gamma\in\Gamma$, $I_{\gamma}^+=I_{\min\{\gamma'\mid\gamma<\gamma'\}}$ for  each
$\gamma\in\Gamma$, $A_{\gamma}=A/I_{\gamma}^+$ for each $\gamma\in\Gamma$, $E_{\gamma}=E/I_{\gamma}^+E$ for  each
$\gamma\in\Gamma$, let $\gr_{\gamma}(A)=I_{\gamma}/I_{\gamma}^+$ for  each $\gamma\in\Gamma$, let
$\gr_{\gamma}(E)=I_{\gamma}E/I_{\gamma}^+E$ for  each $\gamma\in\Gamma$, let $\gr(A)$ be the associated graded ring of
the filtration $\{I_{\gamma}\mid\gamma\in\Gamma\}$, let $\gr(E)$ be the associated graded module of the filtration
$\{I_{\gamma}E\mid\gamma\in\Gamma\}$, let
\[
   \mu_{\gamma}:\gr_{\gamma}(A)\otimes_{A_0}E_0\to\gr_{\gamma}(E)
\]
be the multiplication map for each $n\geq0$, and let
\[
   \mu:\gr(A)\otimes_{A_0}E_0\to\gr(E)
\]
be the direct sum of the $\mu_{\gamma}$.

\begin{thm}\emph{[Local Criterion for t-Flatness]}\label{T:local_criterion}
   Continuing the notation above, $\gr(A)\cong A$ and if $B$ is a Noetherian $A$-algebra, $IB\subseteq\rad(B)$, and $E$
   is a finitely generated $B$-module or if $I$ is nilpotent, then the following are equivalent:
   \begin{itemize}
      \item[(1)] $E$ is t-flat relative to $K$.
      \item[(2)] $\Tor_i^A(N,E)=0$ for all $i>0$ and every combinatorial $A$-module $N$.
      \item[(3)] The canonical surjection $$I\otimes_AE\to IE$$ is an isomorphism.
      \item[(3')] $\Tor_1^A(A_0,E)=0$.
      \item[(4)] $\mu_{\gamma}$ is an isomorphism for all $\gamma\in\Gamma$.
      \item[(4')] $\mu$ is an isomorphism.
      \item[(5)] $E_{\gamma}$ is t-flat relative to $K_{\gamma}^+\cup K$ for all $\gamma\in\Gamma$.
      \item[(6)] The multiplication map $$(I^n/I^{n+1})\otimes_{A_0}E_0\to I^nE/I^{n+1}E $$ is an isomorphism for all
         $n\in\N$.
      \item[(7)] $E/I^{n+1}E$ is t-flat relative to $(n+1)P^+\cup K$ for all $n\in\N$.
   \end{itemize}
   In fact,
   (1)$\Leftrightarrow$(2)$\Rightarrow$(3)$\Leftrightarrow$(3')$\Leftrightarrow$(4)$\Leftrightarrow$(4')$\Leftrightarrow$
   (5)$\Leftrightarrow$(6)$\Leftrightarrow$(7) without any extra assumptions on $I$ or $E$.
\end{thm}

\begin{proof}
   First, consider the underlying group of $\gr(A)$:
   \[
      \gr(A)=\bigoplus_{\gamma\in\Gamma}I_{\gamma}/I_{\gamma}^+\cong\bigoplus_{\gamma\in\Gamma}\left(\bigoplus_{
      \substack{\nu(p)=\gamma\\ p\in P\setminus K}}\Z t^p\right)\cong\bigoplus_{p\in P\setminus K}\Z t^p\cong A.
   \]
   since each ideal $I_{\gamma}$ is $P$-homogeneous. Furthermore, multiplication is given by
   \[
      t^p\cdot t^q=
      \begin{cases}
         t^{p+q}  &\text{if $p+q\notin K$,}\\
         0        &\text{otherwise,}
      \end{cases}
   \]
   in $\gr(A)$ since $\nu(p+q)=\nu(p)+\nu(q)$. Since $\{t^p\mid p\in P\setminus K\}$ is a $\Z$-basis for the free group
   $\gr(A)$, $\gr(A)\cong A$ as rings as well.

   (2)$\Rightarrow$(1): Evident.

   (1)$\Rightarrow$(2): First, we will treat the finitely generated combinatorial $A$-modules. We proceed by induction
   on $i$. Let $i=1$. If $N$ is a finitely generated combinatorial $A$-module, then $N$ has a filtration
   \[
      0=N_0\subset N_1\subset\cdots\subset N_n=N
   \]
   with each $N_l/N_{l-1}\cong\Z[P]/(J_l)$ for some prime ideal $J_l\subset P$ according to Corollary~\ref{C:filtration}.
   Since each $N_l$ is a submodule of $N$, each $N_l$ is annihilated by $(K)$. Since, each $\Z[P]/(J_l)\cong
   N_l/N_{l-1}$, each $\Z[P]/(J_l)$ is
   annihilated by $(K)$ as well. So, each $J_l$ contains $K$. Now we proceed by induction on the length of our
   filtration. If $n=1$, then $N=N_1\cong\Z[P]/(J_1)$ and we are done. If $n>1$, assume $\Tor_1^A(N',E)=0$ for every
   combinatorial $A$-module $N'$ whose filtration has length $n-1$. We have
   \[
      \xymatrix@C=0.5cm{0 \ar[r] & N_{n-1} \ar[rr] && N \ar[rr] && \Z[P]/(J_n) \ar[r] & 0 }.
   \]
   Once we take the tensor product with $E$, we obtain the long exact sequence
   \begin{multline*}
      \xymatrix{\cdots \ar[r] & \Tor_1^A(N_{n-1},E) \ar[r] & \Tor_1^A(N,E)} \\
      \xymatrix{\ar[r] & \Tor_1^A(\Z[P]/(J_n),E) \ar[r] & \cdots}
   \end{multline*}
   Since $E$ is t-flat relative to $K$, $\Tor_1^A(\Z[P]/(J_n),E)=0$. Furthermore,
   \[
      \Tor_1^A(N_{n-1},E)=0
   \]
   by the induction hypothesis. So, we have $\Tor_1^A(N,E)=0$ for every finitely generated combinatorial
   $A$-module $N$. Now suppose $i>1$ and
   \[
      \Tor_{i-1}^A(N',E)=0
   \]
   for every finitely generated combinatorial $A$-module $N'$. Let $J\subseteq P$ be an ideal containing
   $K$. We have
   \[
      \xymatrix@C=0.5cm{0 \ar[r] & (J)/(K) \ar[rr] && A \ar[rr] && \Z[P]/(J) \ar[r] & 0 }.
   \]
   Once we take the tensor product with $E$, we obtain the long exact sequence
   \begin{multline*}
      \xymatrix{\cdots \ar[r] & 0 \ar[r] & \Tor_i^A(\Z[P]/(J),E)} \\
      \xymatrix{\ar[r] & \Tor_{i-1}^A((J)/(K),E) \ar[r] & 0 \ar[r] & \cdots}
   \end{multline*}
   Since $\Tor_{i-1}^A((J)/(K),E)=0$ by the induction hypothesis,
   \[
      \Tor_i^A(\Z[P]/(J),E)=0.
   \]
   Any combinatorial $A$-module is the union of its finitely generated combinatorial submodules. Since right
   exact functors commute with colimits and filtered colimits are exact in module categories, we are done.

   (3)$\Leftrightarrow$(3'): Notice that $\Tor_1^A(A_0,E)$ is the kernel of the surjective map in
   (3).

   (1)$\Rightarrow$(3'): (1) implies (3') follows from the definition of t-flatness.

   (3')$\Rightarrow$(1): Let $S=\{J\mid\Tor_1^A(\Z[P]/(J),E)\neq0\}$. If $S$ is nonempty, then $S$ has a maximal
   element $J$ since the ideals of $P$ satisfy the ascending chain condition. By Lemma~\ref{L:prime}, $J$ is prime. Let
   $q\in P^+\setminus J$. Since $J$ is prime, we have a short exact sequence
   \[
      \xymatrix@C=0.5cm{0 \ar[r] & \Z[P]/(J)\ \ar[rr]^{\cdot t^q} && \Z[P]/(J) \ar[rr] && \Z[P]/(q,J) \ar[r] & 0}
   \]
   This sequence yields the long exact $\Tor$ sequence
   \[
      \xymatrix@C=0.5cm{\cdots \ar[rr] && \Tor_1^A(\Z[P]/(J),E) \ar[rr]^{\cdot t^q} && \Tor_1^A(\Z[P]/(J),E) \ar[r] & 0}
   \]
   Now apply Nakayama's Lemma. We conclude
   \[
      \Tor_1^A(\Z[P]/(J),E)_{\p}=0
   \]
   at every prime $\p$ containing an element of $P^+\setminus J$. If $I$ is nilpotent, then every prime contains every
   element of $P^+\setminus J$. If $E$ is a finitely generated module over some $A$-algebra $B$ such that
   $IB\subseteq\rad(B)$, every maximal ideal in the support of $E$ contains every element of $P^+\setminus J$. In
   either case,
   \[
      \Tor_1^A(\Z[P]/(J),E)_{\m}=0
   \]
   at every maximal ideal $\m\subset A$ in the support of $E$. Therefore,
   \[
      \Tor_1^A(\Z[P]/(J),E)=0
   \]
   and $E$ is t-flat relative to $K$.

   (2)$\Rightarrow$(4):  Mimic the (3)$\Rightarrow$(4) argument in Matsumura.

   (4)$\Leftrightarrow$(4'): Evident.

   (4)$\Rightarrow$(5): Mimic the (4)$\Rightarrow$(5) argument in Matsumura.

   (5)$\Rightarrow$(1): Mimic the (5)$\Rightarrow$(1) argument in Matsumura.

   (3')$\Leftrightarrow$(4)$\Leftrightarrow$(5): Fix $\gamma$ and apply the previous arguments with $K$ replaced by
   $K'=K_{\gamma}^+\cup K$ and $E$ replaced by $E/(K')E$. In this case, $(P^+)/(K')$ is nilpotent. So, (4) and (5) are
   both equivalent to
   \[
      \Tor_1^{\Z[P]/(K')}(A_0,E/(K')E)=0.
   \]
   But, $\Tor_1^{\Z[P]/(K')}(A_0,E/(K')E)$ is the kernel of the canonical map
   \[
      \varphi:(I/I_{\gamma}^+)\otimes_{A_{\gamma}}E_{\gamma}\to E_{\gamma}
   \]
   and
   \[
      (I/I_{\gamma}^+)\otimes_{A_{\gamma}}E_{\gamma}\cong(I/I_{\gamma}^+)\otimes_{A_{\gamma}}A_{\gamma}\otimes_AE\cong
      (I/I_{\gamma}^+)\otimes_AE.
   \]
   So, $\Tor_1^{\Z[P]/(K')}(A_0,E/(K')E)$ is isomorphic to $\Tor_1^A(A_0,E)$ the kernel of
   $$(I/I_{\gamma}^+)\otimes_AE\to E_{\gamma}.$$

   (2)$\Rightarrow$(6): See Matsumura's proof that (3)$\Rightarrow$(4).

   (6)$\Rightarrow$(7): See Matsumura's proof that (4)$\Rightarrow$(5).

   (7)$\Rightarrow$(1): See Matsumura's proof that (5)$\Rightarrow$(1).

   (3')$\Leftrightarrow$(6)$\Leftrightarrow$(7): Mimic our proof that (3')$\Leftrightarrow$(4)$\Leftrightarrow$(5).
\end{proof}

\begin{lemma}\label{L:prime}
   Continuing the above notation, if $J'$ is an ideal of $P$ containing $K$ such that
   \[
      \Tor_1^{\Z[P]/(K)}(\Z[P]/(J'),E)\neq0,
   \]
   then there exists a prime $J$ containing $J'$ such that
   \[
      \Tor_1^{\Z[P]/(K)}(\Z[P]/(J),E)\neq0.
   \]
\end{lemma}

\begin{proof}
   Let $S=\{J''\mid\Tor_1^{\Z[P]/(K)}(\Z[P]/(J''),E)\neq0\}$. Since $S$ is nonempty, $S$ has a maximal element $J$ since
   the ideals of $P$ satisfy the ascending chain condition. Suppose $J$ is an ideal of $P$ containing $J'$ such that
   \[
      \Tor_1^{\Z[P]/(K)}(\Z[P]/(J),E)\neq0
   \]
   and no ideal properly containing $J$ has this property. We will prove $J$ is prime by contradiction. Suppose
   $J$ were not prime. Using Theorem~\ref{T:prime_filtration}, write
   \[
      0=N_0\subset N_1\subset\cdots\subset N_n=\Z[P]/(J)
   \]
   such that each $N_{i+1}/N_i\cong\Z[P]/(J_i)$ for some prime ideal $J_i\in P$. As before, each $J_i$ contains
   $J$. Since $J$ is assumed to not be prime, these containments are proper. We get a series of short exact
   sequences
   \[
      \xymatrix@C=0.5cm{0 \ar[r] & N_i \ar[rr] && N_{i+1} \ar[rr] && \Z[P]/(J_i) \ar[r] & 0}
   \]
   Take the various long exact $\Tor$ sequences to get
   \[
      \xymatrix@C=0.5cm{0 \ar[r] & \Tor_1^{\Z[P]/(K)}(N_i,E) \ar[rr] && \Tor_1^{\Z[P]/(K)}(N_{i+1},E) \ar[r] & 0 }
   \]
   by the maximality of $J$. That is
   \begin{align*}
      \Tor_1^{\Z[P]/(K)}(\Z[P]/(J),E)   &=\Tor_1^{\Z[P]/(K)}(N_n,E)\\
                                    &\cong\Tor_1^{\Z[P]/(K)}(N_0,E)\\
                                    &=0
   \end{align*}
   This contradicts the choice of $J$. So, $J$ must be prime.
\end{proof}

In particular, if $A=\Z[P]/(K)$, an $A$-module $E$ is t-flat (relative to $K$) if and only if $E_{\p}$ is flat over
$A_{\p}$ for every log prime ideal $\p=(J)/(K)$ of $A$.

\begin{prop}\label{P:pushforward}
   Let $P\to Q$ be an inclusion between finitely generated torsion-free monoids, let $K$ be an ideal of $Q$, let
   $K'=K\cap P$, let $N$ be a combinatorial $\Z[Q]/(K)$-module, and let $E$ be a $\Z[P]/(K')$-module. If $E$ has t-flat
   dimension $d$ relative to $K$, then
   \[
      \Tor_i^{\Z[P]/(K')}(N,E)=0\text{ for all }i>d
   \]
\end{prop}

\begin{proof}
   It suffices to prove that if $E$ is t-flat relative to $K$, then $$\Tor_1^{\Z[P]/(K')}(N,E)=0$$ (For $d>0$, apply the
   t-flat case to the $d$th syzygy module of $E$). According to Proposition~\ref{P:inclusion}, $N$ is a direct sum of
   combinatorial $\Z[P]/(K')$-modules. Since tor functors commute with direct sums, $$\Tor_1^{\Z[P]/(K')}(N,E)$$ is the
   direct sum of modules of the form $\Tor_1^{\Z[P]/(K')}(N',E)$ where each $N'$ is a combinatorial $\Z[P]/(K')$-module.
   Each $\Tor_1^{\Z[P]/(K')}(N',E)=0$ by the equivalence of (1) and (2) in Theorem~\ref{T:local_criterion}.
\end{proof}

\begin{prop}\label{P:base_change}
   Let $P\to Q$ be an inclusion between finitely generated torsion-free monoids, let $E$ be a $\Z[P]/(K)$-module, and let
   $K'=K+Q$. If $E$ has t-flat dimension $d$ relative to $K$, then the $\Z[Q]/(K')$-module
   $E'=\Z[Q]/(K')\otimes_{\Z[P]/(K)}E$ has t-flat dimension less than or equal to $d$ relative to $K'$.
\end{prop}

\begin{proof}
   It suffices to prove that if $E$ is t-flat relative to $K$, then $E'$ is t-flat relative to $K'$ (For $d>0$,
   apply the t-flat case to the $d$th syzygy module of $E$). Now suppose $E$ is t-flat.

   Let $N$ be a combinatorial $\Z[Q]/(K')$-module and let $0\to L\to F\to N\to0$ be an exact sequence of
   $\Z[Q]/(K')$-modules with $F$ free. Tensor this exact sequence with $E'$ to get the long exact sequence
   \begin{multline*}
      \xymatrix@C=0.5cm{\cdots \ar[rr] && 0 \ar[rr] && \Tor_1^{\Z[Q]/(K')}(N,E')} \\
      \xymatrix@C=0.5cm{\ar[rr] && L\otimes_{\Z[Q]/(K')}E' \ar[rr] && F\otimes_{\Z[Q]/(K')}E'} \\
      \xymatrix@C=0.5cm{\ar[rr]&& N\otimes_{\Z[Q]/(K')}E' \ar[r] & 0 }.
   \end{multline*}
   Our long exact sequence can also be written as
   \begin{multline*}
      \xymatrix@C=0.5cm{\cdots \ar[rr] && 0 \ar[rr] && \Tor_1^{\Z[Q]/(K')}(N,E')} \\
      \xymatrix@C=0.5cm{\ar[rr] && L\otimes_{\Z[P]/(K)}E' \ar[rr] && F\otimes_{\Z[P]/(K)}E'} \\
      \xymatrix@C=0.5cm{\ar[rr]&& N\otimes_{\Z[P]/(K)}E' \ar[r] & 0 }.
   \end{multline*}
   since, for any $\Z[Q]/(K')$-module $E''$,
   \[
      E''\otimes_{\Z[Q]/(K')}E'=E''\otimes_{\Z[Q]/(K')}\Z[Q]/(K')\otimes_{\Z[P]/(K)}E\cong E''\otimes_{\Z[P]/(K)}E
   \]
   On the other hand, we have the long exact sequence obtained by taking the tensor product of our short exact sequence
   with $E$ over $\Z[P]/(K)$:
   \begin{multline*}
      \xymatrix@C=0.5cm{\cdots \ar[rr] && \Tor_1^{\Z[P]/(K)}(F,E) \ar[rr] && \Tor_1^{\Z[P]/(K)}(N,E)} \\
      \xymatrix@C=0.5cm{\ar[rr] && L\otimes_{\Z[P]/(K)}E' \ar[rr] && F\otimes_{\Z[P]/(K)}E'} \\
      \xymatrix@C=0.5cm{\ar[rr]&& N\otimes_{\Z[P]/(K)}E' \ar[r] & 0 }.
   \end{multline*}
   By Proposition~\ref{P:pushforward} and the fact that tor functors commute with direct sums,
   \[
      \Tor_1^{\Z[P]/(K)}(F,E)=0.
   \]
   Therefore,
   \[
      \Tor_1^{\Z[Q]/(K')}(N,E')\cong\Tor_1^{\Z[P]/(K)}(N,E)=0
   \]
   and $E'$ is t-flat relative to $K'$.
\end{proof}

\begin{defn}
   We say a $\Z[P]/(K)$-module $E$ is \emph{weakly t-flat relative to $K$} if
   \[
      \Tor_1^{\Z[P]/(K)}(\Z[P]/(P^+),E)=0.
   \]
\end{defn}

\begin{lemma}\label{L:weakly}
   Continuing the notation from Theorem~\ref{T:local_criterion}, let $B$ be a Noetherian $A$-algebra, let $E$ be a
   finitely generated $B$-module and let $\widehat{E}$ be the $I$-adic completion of $E$. Then $E$ is weakly t-flat
   relative to $K$ if and only if $\widehat{E}$ is t-flat relative to $K$.
\end{lemma}

\begin{proof}
   First, suppose $E$ is weakly t-flat relative to $K$. Since $B$ is Noetherian, so is its $I$-adic completion
   $\widehat{B}$. Furthermore, $I\widehat{B}\subseteq\rad(\widehat{B})$. Consider the long exact tor sequence
   \[
      \xymatrix@C=0.5cm{\cdots \ar[rr] && 0 \ar[rr] && \Tor_1^A(A_0,E) \ar[rr] && I\otimes_AE \ar[rr] && E_0
         \ar[r] & 0 }.
   \]
   It is an exact sequence of $B$-modules and $\widehat{B}$ is a flat $B$-module. So, after taking the tensor
   product with $\widehat{B}$, the resultant sequence is exact. In particular,
   \[
      \Tor_1^A(A_0,\widehat{E})\cong\Tor_1^A(A_0,E)\otimes_B\widehat{B}=0.
   \]
   We obtain the t-flatness of $\widehat{E}$ by the equivalence of (1) and (3') in
   Theorem~\ref{T:local_criterion}.

   On the other hand, suppose $\Tor_1^A(A_0,\widehat{E})=0$. Let $S$ be the multiplicatively closed set $1+IB$, then
   $\widehat{B}$ is faithfully flat over $S^{-1}B$ and
   \[
      \Tor_1^A(A_0,E)\otimes_B\widehat{B}\cong\Tor_1^A(A_0,\widehat{E})=0.
   \]
   So, $\Tor_1^A(A_0,E)\otimes_BS^{-1}B=0$. But,
   \[
      \Supp(\Tor_1^A(A_0,E))\subseteq\{\p\subset\Z[P]/(K)\mid I\subseteq\p\}.
   \]
   That is, $\Tor_1^A(A_0,E)=0$.
\end{proof}

\begin{lemma}\label{L:upward}
   Let $K\subseteq J$ be proper ideals of $P$ and let $E$ be a $\Z[P]/(K)$-module. If $E$ is weakly t-flat relative to
   $K$, then $E/(J)E$ is weakly t-flat relative to $J$.
\end{lemma}

\begin{proof}
   Continue the notation of Theorem~\ref{T:local_criterion}. By the equivalence of (3) and (3') in
   Theorem~\ref{T:local_criterion}, it suffices to prove that the canonical surjective homomorphism
   \[
      (I/(J))\otimes_{A/(J)}(E/(J)E)\to IE/(J)E
   \]
   is an isomorphism. Since
   \[
      (I/(J))\otimes_AE\cong(I/(J))\otimes_{A/(J)}(A/(J))\otimes_AE\cong(I/(J))\otimes_{A/(J)}(E/(J)E),
   \]
   it suffices to prove that the canonical map $\varphi:(I/(J))\otimes_AE\to E/(J)E$ is injective.
   Take the tensor product of the short exact sequence
   \[
      \xymatrix@C=0.5cm{0 \ar[r] & I/(J) \ar[rr] && A/(J) \ar[rr] && A_0 \ar[r] & 0 }
   \]
   with $E$ over $A$ to see $\ker\varphi\cong\Tor_1^A(A_0,E)$. Since $E$ is weakly t-flat relative to $K$,
   $\ker\varphi=0$.
\end{proof}

\begin{lemma}\label{L:intersection}
   Let $K$, $K_1$ and $K_2$ be proper ideals of $P$, let $K_1\cap K_2\subseteq K$, and let $E$ be a $\Z[P]/(K)$-module.
   If $E/(K_1)E$ is weakly t-flat relative to $K_1$ and $E/(K_2)E$ is weakly t-flat relative to $K_2$, then $E$ is
   weakly t-flat relative to $K$.
\end{lemma}

\begin{proof}
   We continue the notation from Theorem~\ref{T:local_criterion}. It suffices to prove the multiplication map
   \[
      (I^n/I^{n+1})\otimes_{A_0}E_0\to I^nE/I^{n+1}E
   \]
   is an isomorphism for all $n$ by the equivalence of (3) and (7). For any particular $n$, we may assume
   $(n+1)P^+\subseteq K$ when attempting to prove this map is an isomorphism. From now on, assume $(n+1)P^+\subseteq
   K$. By the equivalence of (3) and (3') it suffices to prove $I\otimes E\to E$ is injective. Take the tensor product
   of $E$ with the short exact sequence
   \begin{equation}
      \xymatrix@C=0.5cm{
        0 \ar[r] & I \ar[rr] && I/(K_1)\oplus I/(K_2) \ar[rr] && I/(K_1\cup K_2) \ar[r] & 0 }\label{E:short}
   \end{equation}
   and consider the commutative diagram
   \begin{equation}
      \xymatrix{
         I\otimes E \ar[d]_{f} \ar[r]^-g & (I/(K_1))\otimes E\oplus(I/(K_2))\otimes E  \ar[d]_{h} \ar[r] &
            (I/(K_1\cup K_2))\otimes E \ar[r] & 0 \\
         E \ar[r]^-j & E/(K_1)E\oplus E/(K_2)E & & }\label{E:diagram}
   \end{equation}
   We want to prove the $f$ is injective. Since $E/(K_1)E$ is weakly t-flat relative to $K_1$ and $E/(K_2)E$ is weakly
   t-flat relative to $K_2$, $h$ is injective. Since $h$ is injective, it suffices to prove $g$ is injective, then
   $h\circ g=j\circ f$ is injective and $f$ is injective. We will prove $g$ is injective by induction on $n$. If $n=1$,
   $K=K_1=K_2=P^+$ and we are done. Now suppose $E/(J)E$ is weakly t-flat whenever $nP^+\cup(K_1\cap K_2)\subseteq J$.
   In particular, let $J=K\cup nP^+$. Notice that $I$, $I/(K_1)\oplus I/(K_2)$ and $I/(K_1\cup K_2)$ are combinatorial
   $\Z[P]/(J)$-modules since $I$ is annihilated by $(J)$. Furthermore, the short exact sequence~(\ref{E:short}) above
   is an exact sequence of $\Z[P]/(J)$-modules. In fact, the top line of our commutative diagram~(\ref{E:diagram}) can
   be obtained by taking the tensor product over $\Z[P]/(J)$ of our short exact sequence~(\ref{E:short}) above with
   $E/(J)E$. So,
   \[
      \Tor_1^A(I/(K_1\cup K_2),E)=\Tor_1^{\Z[P]/(J)}(I/(K_1\cup K_2),E/(J)E)
   \]
   Since $E/(J)E$ is weakly t-flat relative to $J$ and $I/(K_1\cup K_2)$ is a combinatorial
   $\Z[P]/(J)$-module, we have $\Tor_1^{\Z[P]/(J)}(I/(K_1\cup K_2),E/(J)E)=0$. Therefore, $g$ is an injection.
\end{proof}

\begin{prop}\label{P:ideal}
   Let $A$ be a Noetherian $\Z[P]$-algebra and let $E$ be a finitely generated $A$-module. There exists an ideal
   $K\subseteq P^+$ such that the module $E/(J)E$ is weakly t-flat relative to an ideal $J$ if and only if $K\subseteq
   J$.
\end{prop}

\begin{proof}
   Let $S=\{J\mid E/(J)E\text{ is weakly t-flat relative to }J\}$, let $K=\bigcap_{J\in S}J$, and let $I=(P^+)/(K)$. By
   Lemma~\ref{L:upward}, it suffices to prove $E/(K)E$ is weakly t-flat relative to $K$, for then $S=\{J\mid K\subseteq
   J\}$. By the equivalence of (3') and (7) in Theorem~\ref{T:local_criterion}, to prove $E/(K)E$ is weakly t-flat
   relative to $K$, it suffices to prove $E/I^{n+1}E$ is t-flat relative to $(n+1)P^+\cup K$ for all $n\in\N$. If $J\in
   S$, then $E/(I^{n+1}+(J))E$ is t-flat relative to $(n+1)P^+\cup J$ for all $n\in\N$ by the same equivalence. That
   is, $(n+1)P^+\cup J\in S$ for all $n\in\N$. Furthermore, $(n+1)P^+\cup K=\bigcap_{J\in S}((n+1)P^+\cup J)$. Since
   there are only finitely many ideals of $P$ containing $(n+1)P^+$, this intersection is finite and $E/I^{n+1}E$ is
   t-flat relative to $(n+1)P^+\cup K$ by Lemma~\ref{L:intersection}.
\end{proof}

\subsection{Openness of t-Flat Loci}

\begin{defn}\cite[Definition~2.4]{aO97}
   Let $(X,\M)$ be a locally Noetherian coherent log scheme. We say a sheaf of ideals
   $\K\subseteq\M$ is \emph{coherent} if, locally on $X$, there exists a fine chart $P_X\to\M$ and an ideal
   $K\subseteq P$ such that $\K$ is the ideal sheaf generated by the image of $K$.
\end{defn}

\begin{prop}\emph{\cite[Proposition~2.6]{aO97}}\label{P:coherent_ideal}
   Let $(X,\M)$ be a locally Noetherian fine log scheme, let $\K\subseteq\M$ be a coherent
   sheaf of ideals, let $\beta:P_X\to\M$ be a fine chart, and let $K=\beta^{-1}(\K)$. Then $\K$ is the
   ideal sheaf generated by the image of $K$.
\end{prop}

\begin{thm}\label{T:open}
   Let $(X,\M)$ be a locally Noetherian toroidal log scheme, let $\K\subseteq\M$ be a coherent
   ideal sheaf, and let $\F$ be a coherent sheaf of $\O_X$-modules annihilated by $(\K)$. Then
   \[
      \{x\in X\mid\F\text{ is $\M$-flat relative to }\K\text{ at }x\}
   \]
   is open.
\end{thm}

\begin{proof}
   The question is local. Fix $x\in X$. By Proposition~\ref{P:canonical_chart}, we may assume
   $X=\Spec(P\xrightarrow{\beta}A)$ where $P=\overline{\M}_x$, $A$ is a Noetherian ring and $x=\p$. In
   particular, $\beta$ is local at $x$ and we have an associated ring homomorphism $\widetilde{\beta}:\Z[P]\to
   A$. By the previous proposition, we may also assume $\K$ is generated by the image of an ideal $K\subseteq
   P$. Furthermore, we may assume $\F$ is the sheaf associated to some finitely generated $A$-module $E$.

   By Corollary~\ref{C:comparison}, $\F$ is $\M$-flat relative to $\K$ at $x\in X$ if and only if $\F_x$ is
   t-flat relative to $K$. So, $\F$ is not $\M$-flat relative to $\K$ at $x\in X$ if and only if
   $\widetilde{\beta}^{-1}(\p)$ is in the support of $\Tor_1^{\Z[P]/(K)}(\Z[P]/(J),E)$ for some ideal $J\subset P$
   containing $K$.  Since $\Z[P]/(J)$ is a finitely generated combinatorial $\Z[P]$-module, in order to check
   \[
      \Tor_1^{\Z[P]/(K)}(\Z[P]/(J),E)=0\quad\forall J\subset P\text{ containing }K,
   \]
    by Lemma~\ref{L:prime}, it suffices to check only the prime ideals $J\subset P$ containing $K$. Since
   $P$ is finitely generated, $P$ has only finitely many primes. Therefore,
   \[
      \{\Tor_1^{\Z[P]/(K)}(\Z[P]/(J),E)\mid K\subseteq J\subset P,\,J\text{ prime}\}
   \]
   is a finite set of finitely generated modules. So, the union of the supports of these modules is closed. But
   $x$ is in this union if and only if $\F$ is not $\M$-flat relative to $\K$ at $x$. So,
   \[
      \{x\in X\mid\F\text{ is $\M$-flat relative to }\K\text{ at }x\}
   \]
   is open.
\end{proof}

\begin{prop}\label{P:finite}
   Let $P$ be a finitely generated torsion-free monoid, let $A$ be a Noetherian ring, let $\beta:P\to A$ be a
   monoid homomorphism with respect to multiplication, let $E$ be a finitely generated $A$-module, and let
   $K_E:\Spec A\to\{\text{ideals of $P$}\}$ be the function that takes a prime $\p$ to the ideal $K_E(\p)$ of $P$
   such that $E_{\p}$ is weakly t-flat relative to $J$ if and only if $K_E(\p)\subseteq J$. Then, the image of
   $K_E$ is finite.
\end{prop}

\begin{proof}
   We will prove this by Noetherian induction. It suffices to prove that $K_E$ is constant on some nonempty open
   subset of $\Spec A$. Let $\{X_i\}_{i=1}^n$ be the set of irreducible components of $\Spec A$, let
   $X=X_1\setminus\bigcup_{i=2}^nX_i$ and let $\eta$ be the generic point of $X$. Let $K$ be the ideal of $P$ such
   that $E_{\eta}$ is weakly t-flat relative to $J$ if and only if $K\subseteq J$. Now consider $X_K=\{x\in X\mid
   K_E(x)=K\}$. Since $\eta$ is in $X_K$, $X_K$ is nonempty. For each point $x$ on $X$, let $U(x)=\{x'\in X\mid
   E_{x'}\text{ is weakly t-flat relative to }K_E(x)\}$. Note that $U(x)=\{x'\in X\mid K_E(x')\subseteq
   K_E(x)\}$. By Theorem~\ref{T:open}, $U(x)$ is open for all $x$. Since each $U(x)$ is open, $\eta$ is an
   element of $U(x)$ for all $x$ in $X$. That is, $K_E(\eta)\subseteq K_E(x)$ for all $x$ in $X$. So,
   $X_K=U(\eta)$ and $X_K$ is open.
\end{proof}

\begin{xmpl}
   Let $P=\N^2$, let $A=\C[x,y]$, let $\beta:P\to A$ be given by $(n,m)\mapsto x^ny^m$, and let $E=A/(x-y)\bigoplus
   A/(y^3)$. Here
   \[
      K_E(\p)=
      \begin{cases}
         P^+,     &\text{if $\p=(x,y)$;} \\
         (0,3)+P  &\text{if $\p=(x-\alpha,y)$ with $\alpha\neq0$ or $\p=(y)$;} \\
         \varnothing &\text{otherwise.}
      \end{cases}
   \]
\end{xmpl}

\section{Log Regularity}\label{Ch:regularity}

This section generalizes much of Kato~\cite{kK94} by relaxing his condition (S) to condition $(*)$ below. By doing so,
we work with toriodal log schemes rather than fine saturated log schemes. Most of Kato's methods go through with only
minor modifications.

\subsection{Definition of Toric Singularity}

In this chapter, we will mainly consider log schemes which satisfy the following condition $(*)$.

$(*)$: $(X,\M)$ is toriodal and its underlying scheme $X$ is locally Noetherian.

\begin{defn}
   Let $(X,\M)$ be a log scheme satisfying condition $(*)$. We say $(X,\M)$ is
   \emph{logarithmically regular} at $x$, or $(X,\M)$ has (at worst) a \emph{toric singularity} at $x$,
   if the following two conditions are satisfied.
   \begin{itemize}
      \item[(i)] $\O_{X,x}/(\M_x^+)$ is a regular local ring.
      \item[(ii)] $\dim(\O_{X,x})=\dim(\O_{X,x}/(\M_x^+))+\rank(\overline{\M}_x^{gp})$.
   \end{itemize}
   We say $(X,\M)$ is \emph{log regular} if $(X,\M)$ is log regular at each point $x \in X$.
\end{defn}

\begin{lemma}
   \emph{\cite[Lemma~(2.3)]{kK94}} Let $(X,\M)$ be a log scheme satisfying $(*)$ and let $x\in
   X$. Then
   \[
      \dim (\O_{X,x})\leq\dim(\O_{X,x}/(\M_x^+))+\rank(\overline{\M}_x^{gp}).
   \]
\end{lemma}

\begin{proof}
   See Kato's proof.
\end{proof}

\subsection{Completed Toric Singularities}

In this subsection, continue to let $(X,\M)$ be a log scheme satisfying $(*)$.

\begin{lemma}\label{L:injective}
   \emph{\cite[Lemma~(3.5)]{kK94}} Let $R$ be a ring, let $\pi$ be a nonzero-divisor of R, let $P$ and $Q$ be
   affine semigroups, and let $P\to Q$ be an injective homomorphism. Let $\theta$ be an element of $R[[P]]$ such
   that $\theta\equiv\pi\bmod{(P^+)}$. Then $R[[P]]/(\theta)\to R[[Q]]/(\theta)$ is injective.
\end{lemma}

\begin{proof}
   See Kato's proof.
\end{proof}

\begin{lemma}\label{L:domain}
   \emph{\cite[Lemma~(3.4)]{kK94}} Let $R$ be a ring, let $\pi$ be a nonzero-divisor of $R$ such that $R/(\pi)$
   is an integral domain, let $P$ be an affine semigroup. Let $\theta$ be an element of
   $R[[P]]$ such that $\theta\equiv\pi\bmod{(P^+)}$. Then $R[[P]]/(\theta)$ is an integral domain.
\end{lemma}

\begin{proof}
   See Kato's proof.
\end{proof}

\begin{thm}\label{T:iso1}
   \emph{\cite[Theorem~(3.2)]{kK94}} Let $x\in X$. Assume $\O_{X,x}/(\M_x^+)$ is regular, let
   $P=\overline{\M}_x$, let $\varphi$ be a section of $\M_x\to P$ as in Corollary~\ref{C:chart}, and let
   $f_1,\ldots,f_d\in\O_{X,x}$ such that $(f_i\bmod{\M_x^+})_{1\leq i\leq d}$ is a regular system of parameters
   of $\O_{X,x}/(\M_x^+)$.
   \begin{enumerate}
      \item If $\O_{X,x}$ contains a field, let $k$ be a subfield of $\widehat{\O}_{X,x}$ such that
         $k\cong\widehat{\O}_{X,x}/\widehat{\m}_x$. Then $(X,\M)$ is log regular at $x$ if and
         only if the surjective homomorphism
         \[
            \psi: k[[P]][[t_1,\ldots,t_d]]\to\widehat{\O}_{X,x};\,t_i\mapsto f_i
         \]
         is an isomorphism.
      \item If $\O_{X,x}$ does not contain a field, let $R$ be a complete discrete valuation ring in which
         $p=\operatorname{char}(\O_{X,x}/\m_x)$ is a prime element and fix a homomorphism
         $R\to\widehat{\O}_{X,x}$ which induces $R/pR\xrightarrow{\cong}\widehat{\O}_{X,x}/\widehat{\m}_x$. Then
         $(X,\M)$ is log regular at $x$ if and only if the kernel of the surjective homomorphism
         \[
            \psi:R[[P]][[t_1,\ldots,t_d]]\to\widehat{\O}_{X,x};\,t_i\mapsto f_i
         \]
         is generated by an element $\theta$ such that
         \[
            \theta\equiv p\bmod{(P^+,t_1,\ldots,t_d)}
         \]
    \end{enumerate}
\end{thm}

\begin{proof}
   See Kato's proof.
\end{proof}

\begin{cor}\label{C:iso2}
   \emph{\cite[Theorem~(3.1)]{kK94}}
   \begin{enumerate}
      \item $(X,\M)$ is log regular at $x$ if and only if there exists a complete regular local
         ring $R$, an affine semigroup $P$, and an isomorphism
         \[
            R[[P]]/(\theta)\xrightarrow{\cong}\widehat{\O}_{X,x}
         \]
         with $\theta\in R[[P]]$ satisfying the following conditions
         (i) and (ii).
         \begin{enumerate}
            \item[(i)] The constant term of $\theta$ belongs to $\m_R\setminus\m_R^2$.
            \item[(ii)] The inverse image of $\M$ on $\Spec(\widehat{\O}_{X,x})$ is induced by the map
               $P\to\widehat{\O}_{X,x}$.
         \end{enumerate}
      \item Assume $\O_{X,x}$ contains a field. Then, $(X,\M)$ is log regular at $x$ if and only
         if there exists a field $k$, an affine semigroup $P$, and an isomorphism
         \[
            k[[P]][[t_1,\ldots,t_d]]\xrightarrow{\cong}\widehat{\O}_{X,x}
         \]
         for some $d\geq 0$ satisfying the condition (ii) in (1).
   \end{enumerate}
\end{cor}

\begin{prop}\label{P:t-flat_regular}
   $(X,\M)$ is log regular at $x$ if and only if $\O_{X,x}/(\M_x^+)$ is a regular local ring and
   $\O_{X,x}$ is t-flat.
\end{prop}

\begin{proof}
   Since $\widehat{\O}_{X,x}$ is faithfully flat over $\O_{X,x}$, $\O_{X,x}$ is t-flat if and only if
   $\widehat{\O}_{X,x}$ is t-flat. By Theorem~\ref{T:iso1}, $(X,\M)$ is log regular at $x$ if and
   only if $\Spec(\widehat{\O}_{X,x})$ is log regular at $x$. So, we may assume $\O_{X,x}=\widehat{\O}_{X,x}$. If
   $\O_{X,x}$ is t-flat, then the equivalence of (1) and (4') in the Theorem~\ref{T:local_criterion} shows
   $\dim(\O_{X,x})=\dim(\O_{X,x}/(\M_x^+))+\rank(\overline{\M}_x^{gp})$. On the other hand, if
   $(X,\M)$ is log regular at $x$, then Theorem~\ref{T:iso1} and the equivalence of (1) and (4')
   in the Theorem~\ref{T:local_criterion} show $\O_{X,x}$ is t-flat.
\end{proof}

\begin{thm} \label{T:flat}
   \emph{\cite[Theorem~(6.2)]{kK94}} Let $P=\overline{\M}_x$ be an affine semigroup and take a
   section $\varphi:P\to\M_x$ of $\M_x\to P$ as in Corollary~\ref{C:chart}. Assume $\O_{X,x}$ contains a
   field $k$. Then $(X,\M)$ is log regular at $x$ if and only if $\O_{X,x}/(\M_x^+)$ is a regular
   local ring and the map $k[P]\to\O_{X,x}$ induced by $\varphi$ is flat.
\end{thm}

\begin{proof}
   This follows immediately from Theorem~\ref{P:t-flat_regular}, Theorem~\ref{T:local_criterion}, and the local
   criterion for flatness (see EGA~III~\cite[$0_\mathrm{{III}}$~(10.2.2)]{GD61}).
\end{proof}

\subsection{Some Properties of Toric Singularities}

Kato reminds us that if $(X,\M)$ is a log regular scheme satisfying his condition (S) then its
underlying scheme is Cohen-Macaulay and normal, see Hochster~\cite{mH72}. Let $k$ be a field. If the canonical log
structure on $\Spec k[P]$ satisfies $(*)$, it need not be Cohen-Macaulay nor normal. Consider the monoids
$$\langle(4,0),(3,1),(1,3),(0,4)\rangle\subset\N^2$$ and $\langle2,3\rangle\subset\N$. Information on when affine
semigroup rings are Cohen-Macaulay and its dependence on the characteristic, can be found in \cite{TH86,TH88}.
Information on the local cohomology modules and dualizing complexes of affine semigroup rings is found in \cite{mI88}
and \cite{SS90}.

\begin{thm}
   \emph{\cite[Theorem~(4.2)]{kK94}} Let $(X,\M)$ and $(Y,\mathscr{N})$ be log regular
   schemes. Let $f:(X,\M)\to(Y,\mathscr{N})$ be a morphism whose underlying map of schemes
   is a closed immersion and assume $f^*\mathscr{N}\cong\M$. Then the underlying map of schemes of $f$ is a regular
   immersion.
\end{thm}

\begin{proof}
   See Kato's proof.
\end{proof}

\subsection{Localization}

For the duration of this subsection, let $(X,\M)$ be a log scheme satisfying condition $(*)$.

\begin{lemma} \label{L:tech}
   Let $A$ be a local ring, let $P$ be a 1-dimensional affine semigroup, let $\beta:P\to A$ be a monoid
   homomorphism with respect to multiplication, and let $\p\subset A$ be a prime ideal such that
   $\beta^{-1}(\p)=\emptyset$. If $\Spec(P\xrightarrow{\beta}A)$ is log regular, then $A_{\p}$ is a regular local
   ring.
\end{lemma}

\begin{proof}
   By Theorem~\ref{T:iso1} $\widehat{A}=R[[P]]/(\theta)$ for some complete regular local ring $R$ and some
   $\theta$ whose constant term is contained in $\m_{R}\setminus\m_{R}^2$. Note that $\N$ is the saturation of
   $P$ and there exists an $n \in\N$ such that $m\in P$ whenever $m\geq n$. Let $F$ be the image of $P$ in $A$ and
   consider the following diagram:
   \[
      \xymatrix{
        A \ar[d] \ar[r] & R[[P]]/(\theta) \ar[d] \ar[r] & R[[\N]]/(\theta) \ar[d] \\
        A_F \ar[d] \ar[r] & (R[[P]]/(\theta))_F \ar[d] \ar[r]^= & R((\N))/(\theta) \\
        A_{\p} \ar[r] & (A\setminus\p)^{-1}R[[P]]/(\theta) &}
   \]
   $A_{\p}$ is faithfully flat under $(A\setminus\p)^{-1}R[[P]]/(\theta)$ and the latter is regular since it is
   a localization of $R((\N))/(\theta)$. Hence $A_{\p}$ is regular by faithfully flat descent, see
   Matsumura~\cite[Theorem~23.7]{hM86}.
\end{proof}

\begin{thm} \label{T:mod}
   \emph{\cite[Proposition~(7.2)]{kK94}} Let $x\in X$ and suppose $(X,\M)$ is log regular at $x$.
   Let $\p$ be a prime ideal of $\M_x$ and endow $X'=\Spec(\O_{X,x}/(\p))$ with the log structure $\M'$ associated
   to $\M_x\setminus\p\to\O_{X,x}/(\p)$. $(X',\M')$ is log regular at $x \in X'$.
\end{thm}

\begin{proof}
   Follow Kato's induction proof, using Lemma~\ref{L:tech} in the one-dimensional case.
\end{proof}

\begin{cor} \label{C:prime}
   \emph{\cite[Corollary~(7.3)]{kK94}} With the notation as in the previous theorem, $(\p)$ is a prime ideal of
   height $\dim(\M_x)_{\p}$.
\end{cor}

\begin{thm} \label{T:gen}
   \emph{\cite[Proposition~(7.1)]{kK94}} Let $x\in X$, and assume that $(X,\M)$ is log regular at
   $x$. Then for any $y\in X$ such that $x\in\overline{\{y\}}$, $(X,\M)$ is log regular at $y$.
\end{thm}

\begin{proof}
   Follow Kato's proof and use Lemma~\ref{L:tech} at the last step.
\end{proof}

\subsection{Log Smooth Morphisms}

Let $(X,\M_X)$ and $(Y,\M_Y)$ be log schemes satisfying $(*)$, and let
$f:(X,\M_X)\to(Y,\M_Y)$ be a morphism. Then, the following two conditions (i) and (ii)
are equivalent. We say $f$ is log smooth if $f$ satisfies these conditions. This equivalence was shown in
Kato~\cite{kK89} for log structures on \'{e}tale sites, and the proof there works for the present situation (log
structures on Zariski sites). See Kato~\cite[Section~3.1.5]{kK89a}. So, we omit the proof here.
\begin{itemize}
   \item[(i)] Assume we are given a commutative diagram of log schemes
      \[
         \xymatrix{
           (T,\M_T) \ar[d]_i \ar[r]^g & (X,\M_X) \ar[d]^f \\
           (T',\M_{T'}) \ar[r]^{g'} & (Y,\M_Y)}
      \]
      such that $(T,\M_T)$ and $(T',\M_{T'})$ satisfy $(*)$, $i:T\to T'$ is a closed immersion, $T$ is defined in $T'$
      by a nilpotent ideal of $\O_{T'}$ and $i^*\M_{T'}\to\M_T$ is an isomorphism. Then, locally on $T'$ there is a
      morphism $h:(T',\M_{T'})\to(X,\M_X)$ such that $h\circ i=g$ and $f\circ h=g'$. Furthermore, the underlying
      morphism of schemes $X\to Y$ is locally of finite type. \item[(ii)] \'{E}tale locally on $X$ and $Y$, there exist
      finitely generated, torsion-free monoids $P$ and $Q$, an injective homomorphism $h:P\to Q$ such that the order of
      the torsion part of $Q^{gp}/ h^{gp}(P^{gp})$ is invertible on $X$, and a commutative diagram of log schemes of
      the form
      \[
         \xymatrix{
           (X,\M_X) \ar[d]_{f} \ar[r] & \Spec(Q\hookrightarrow\Z[Q]) \ar[d]_{\scriptstyle{\text{induced by }}h} \\
           (Y,\M_Y) \ar[r] & \Spec(P\hookrightarrow\Z[P])}
      \]
      such that the inverse image of $P$ on $Y$ is $\M_Y$, the inverse image of $Q$ on $X$ is $\M_X$ and the
      induced morphism of the underlying schemes
      \[
         X\to Y\times_{\Spec\Z[P]}\Spec\Z[Q]
      \]
      is smooth (in the classical sense).
\end{itemize}

\begin{thm} \label{T:smooth}
   \emph{\cite[Theorem~(8.2)]{kK94}} Let $f:(X,\M)\to(Y,\mathscr{N})$ be a log smooth
   morphism between log schemes satisfying $(*)$, and assume $(Y,\mathscr{N})$ is log regular. Then
   $(X,\M)$ is log regular.
\end{thm}

\begin{proof}
   See Kato's proof.
\end{proof}

\begin{thm}
   \emph{\cite[Proposition~(8.3)]{kK94}} Let $k$ be a field and let $(X,\M)$ be a log scheme
   satisfying $(*)$ such that the underlying scheme $X$ is a $k$-scheme which is locally of finite type. Then:
   \begin{enumerate}
      \item If $(X,\M)$ is log smooth over $\Spec(\{0\}\to k)$, then $(X,\M)$
         is log regular.
      \item The converse of (1) is true if $k$ is perfect.
   \end{enumerate}
\end{thm}

\begin{proof}
   See Kato's proof.
\end{proof}

For completeness, we take note of the following proposition of Kato:

\begin{prop}
   \emph{\cite[Proposition~(12.2)]{kK94}} Let $(A,\m)$ be a Noetherian local ring, let $P$ be a sharp finitely
   generated torsion-free monoid, let $\beta:P\to A$ be local, and suppose $(X,\M)=\Spec(P\xrightarrow{\beta}A)$ is log
   regular. If $B=\Z[P^{sat}]\otimes_{\Z[P]}A$ and $(Y,\mathcal{N})=\Spec(P^{sat}\to B)$, then:
   \begin{itemize}
      \item[(1)] $(Y,\mathcal{N})$ is log regular.
      \item[(2)] $B$ is the normalization of $A$ and it is a local ring.
      \item[(3)] If $A$ is already normal, then $P$ is already saturated.
   \end{itemize}
\end{prop}

\providecommand{\bysame}{\leavevmode\hbox to3em{\hrulefill}\thinspace}
\providecommand{\MR}{\relax\ifhmode\unskip\space\fi MR }
\providecommand{\MRhref}[2]{%
  \href{http://www.ams.org/mathscinet-getitem?mr=#1}{#2}
}
\providecommand{\href}[2]{#2}


\begin{thebibliography}{10}\addcontentsline{toc}{section}{References}

\bibitem{AMRT75}
A.~Ash, D.~Mumford, M.~Rapoport, and Y.~Tai, \emph{Smooth compactification of
  locally symmetric varieties}, Math. Sci. Press, Brookline, Mass., 1975, Lie
  Groups: History, Frontiers and Applications, Vol. IV. \MR{56:15642}

\bibitem{vD78}
V.~I. Danilov, \emph{The geometry of toric varieties}, Uspekhi Mat. Nauk
  \textbf{33} (1978), no.~2(200), 85--134, 247, English translation: Russian
  Math. Surveys {\bf 33} (1978), no. 2, 97 - 154. \MR{80g:14001}

\bibitem{pD88}
P.~Deligne, \emph{Lettre \`a {L}. {I}llusie}, 1/6/88, 1988.

\bibitem{mD70}
M.~Demazure, \emph{Sous-groupes alg\'ebriques de rang maximum du groupe de
  {C}remona}, Ann. Sci. \'Ecole Norm. Sup. (4) \textbf{3} (1970), 507--588.
  \MR{44:1672}

\bibitem{dE95}
D.~Eisenbud, \emph{Commutative algebra}, Graduate Texts in Mathematics, vol.
  150, Springer-Verlag, New York, 1995, With a view toward algebraic geometry.
  \MR{97a:13001}

\bibitem{gF90}
G.~Faltings, \emph{${F}$-isocrystals on open varieties: results and
  conjectures}, The Grothendieck Festschrift, Vol. II, Progr. Math., vol.~87,
  Birkh\"auser Boston, Boston, MA, 1990, pp.~219--248. \MR{92f:14015}

\bibitem{rG84}
R.~Gilmer, \emph{Commutative semigroup rings}, Chicago Lectures in Mathematics,
  University of Chicago Press, Chicago, Ill., 1984. \MR{85e:20058}

\bibitem{GD61}
A.~Grothendieck and J.~Dieudonn\'e, \emph{\'{E}l\'ements de g\'eom\'etrie
  alg\'ebrique. {I}{I}{I}. \'{E}tude cohomologique des faisceaux coh\'erents.
  {I}}, Inst. Hautes \'Etudes Sci. Publ. Math. No. \textbf{11} (1961), 167.
  \MR{29:1209}

\bibitem{mH72}
M.~Hochster, \emph{Rings of invariants of tori, {C}ohen-{M}acaulay rings
  generated by monomials, and polytopes}, Ann. of Math. (2) \textbf{96} (1972),
  318--337. \MR{46:3511}

\bibitem{oH88}
O.~Hyodo, \emph{A note on the $p$-adic \'etale cohomology in the semi-stable
  reduction case}, Invent. Math. \textbf{91} (1988), no.~3, 543--557.
  \MR{89h:14017}

\bibitem{oH91}
\bysame, \emph{On the de {R}ham-{W}itt complex attached to a semi-stable
  family}, Compositio Math. \textbf{78} (1991), no.~3, 241--260. \MR{93c:14022}

\bibitem{HK94}
O.~Hyodo and K.~Kato, \emph{Semi-stable reduction and crystalline cohomology
  with logarithmic poles}, Ast\'erisque (1994), no.~223, 221--268, P\'eriodes
  $p$-adiques (Bures-sur-Yvette, 1988). \MR{95k:14034}

\bibitem{lI94}
L.~Illusie, \emph{Logarithmic spaces (according to {K}. {K}ato)}, Barsotti
  Symposium in Algebraic Geometry (Abano Terme, 1991), Perspect. Math.,
  vol.~15, Academic Press, San Diego, CA, 1994, pp.~183--203. \MR{95j:14023}

\bibitem{mI88}
M.-{N}. Ishida, \emph{The local cohomology groups of an affine semigroup ring},
  Algebraic geometry and commutative algebra, Vol.\ I, Kinokuniya, Tokyo, 1988,
  pp.~141--153. \MR{90a:13029}

\bibitem{kK89a}
K.~Kato, \emph{Logarithmic degeneration and {D}ieudonn\'e theory}, Preprint,
  1989.

\bibitem{kK89}
\bysame, \emph{Logarithmic structures of {F}ontaine-{I}llusie}, Algebraic
  analysis, geometry, and number theory (Baltimore, MD, 1988), Johns Hopkins
  Univ. Press, Baltimore, MD, 1989, pp.~191--224. \MR{99b:14020}

\bibitem{kK94b}
\bysame, \emph{Semi-stable reduction and $p$-adic \'etale cohomology},
  Ast\'erisque (1994), no.~223, 269--293, P\'eriodes $p$-adiques
  (Bures-sur-Yvette, 1988). \MR{95i:14020}

\bibitem{kK94}
\bysame, \emph{Toric singularities}, Amer. J. Math. \textbf{116} (1994), no.~5,
  1073--1099. \MR{95g:14056}

\bibitem{KKMS73}
G.~Kempf, F.~Knudsen, D.~Mumford, and B.~Saint-Donat, \emph{Toroidal
  embeddings. {I}}, Springer-Verlag, Berlin, 1973, Lecture Notes in
  Mathematics, Vol. 339. \MR{49:299}

\bibitem{pL01}
P.~Lorenzon, \emph{Indexed algebras associated to a log structure and a theorem
  of $p$-descent on log schemes}, Manuscripta Math. \textbf{101} (2000), no.~3,
  271--299. \MR{2001b:14029}

\bibitem{hM86}
H.~Matsumura, \emph{Commutative ring theory}, Cambridge Studies in Advanced
  Mathematics, vol.~8, Cambridge University Press, Cambridge, 1986, Translated
  from the Japanese by M. Reid. \MR{88h:13001}

\bibitem{MO75}
K.~Miyake and T.~Oda, \emph{Almost homogeneous algebraic varieties under
  algebraic torus action}, Manifolds---Tokyo 1973 (Proc. Internat. Conf.,
  Tokyo, 1973) (Tokyo), Univ. Tokyo Press, 1975, pp.~373--381. \MR{52:406}

\bibitem{wN98}
W.~Niziol, \emph{Toric singularities: Log-blow-ups and global resolutions},
  Preprint, www.math.utah.edu/~niziol, 1998.

\bibitem{tO78}
T.~Oda, \emph{Torus embeddings and applications}, Tata Institute of Fundamental
  Research Lectures on and Physics, vol.~57, Tata Institute of Fundamental
  Research, Bombay, 1978, Based on joint work with Katsuya Miyake.
  \MR{81e:14001}

\bibitem{tO88}
\bysame, \emph{Convex bodies and algebraic geometry}, Ergebnisse der Mathematik
  und ihrer Grenzgebiete (3) [Results in Mathematics and Related Areas (3)],
  vol.~15, Springer-Verlag, Berlin, 1988, An introduction to the theory of
  toric varieties, Translated from the Japanese. \MR{88m:14038}

\bibitem{tO91}
\bysame, \emph{Geometry of toric varieties}, Proceedings of the Hyderabad
  Conference on Algebraic Groups (Hyderabad, 1989) (Madras), Manoj Prakashan,
  1991, pp.~407--440. \MR{92m:14068}

\bibitem{aO97}
A.~Ogus, \emph{Logarithmic {D}e {R}ham cohomology}, Preprint, 1997.

\bibitem{GSR99b}
J.~C. Rosales and P.~A. Garc{\'\i}a-S{\'a}nchez, \emph{Finitely generated
  commutative monoids}, Nova Science Publishers Inc., Commack, NY, 1999.
  \MR{2000d:20074}

\bibitem{iS73}
I.~Satake, \emph{On the arithmetic of tube domains (blowing-up of the point at
  infinity)}, Bull. Amer. Math. Soc. \textbf{79} (1973), 1076--1094.
  \MR{48:8861}

\bibitem{SS90}
U.~Sch{\"a}fer and P.~Schenzel, \emph{Dualizing complexes of affine semigroup
  rings}, Trans. Amer. Math. Soc. \textbf{322} (1990), no.~2, 561--582.
  \MR{92a:13012}

\bibitem{jS76}
J.~H.~M. Steenbrink, \emph{Limits of {H}odge structures}, Invent. Math.
  \textbf{31} (1975/76), no.~3, 229--257. \MR{55:2894}

\bibitem{hS74}
H.~Sumihiro, \emph{Equivariant completion}, J. Math. Kyoto Univ. \textbf{14}
  (1974), 1--28. \MR{49:2732}

\bibitem{hS75}
\bysame, \emph{Equivariant completion. {I}{I}}, J. Math. Kyoto Univ.
  \textbf{15} (1975), no.~3, 573--605. \MR{52:8137}

\bibitem{hT02}
H.~Thompson, \emph{On toric log schemes}, Ph. D. dissertation, 2002.

\bibitem{TH86}
N.~V. Trung and L.~T. Hoa, \emph{Affine semigroups and {C}ohen-{M}acaulay rings
  generated by monomials}, Trans. Amer. Math. Soc. \textbf{298} (1986), no.~1,
  145--167. \MR{87j:13032}

\bibitem{TH88}
\bysame, \emph{Corrigendum to: ``{A}ffine semigroups and {C}ohen-{M}acaulay
  rings generated by monomials''}, Trans. Amer. Math. Soc. \textbf{305} (1988),
  no.~2, 857. \MR{88m:13008}

\end{thebibliography}
\end{document}